\documentclass[Threeside,11pt,tikz,border=3mm]{article}
\usepackage{etoolbox}
\makeatletter
\patchcmd{\thebibliography}
{\list}
{\small          
	\setstretch{1.2}
	\list}
{}{}
\makeatother
\usepackage{setspace} 
\usepackage{tikz}
\usepackage{tikz-cd}
\usetikzlibrary{arrows.meta}
\usepackage{subcaption}
\usepackage [latin1]{inputenc}
\usepackage[english]{babel}
\usepackage{amsmath,amscd}
\usepackage{amsthm}
\usepackage{algorithm}
\usepackage{amsfonts}
\usepackage{amssymb}
\usepackage{mathrsfs}
\usepackage{graphicx}
\usepackage{colortbl,dcolumn}
\usepackage{marvosym}
\usepackage{ifsym}
\usepackage{paralist}
\usepackage{xcolor}
\usepackage{array,color}
\usepackage{geometry}
\usepackage[colorlinks,
citecolor=blue,
linkcolor=blue,
backref=page
allcolors=black]{hyperref} 
\usetikzlibrary{shapes.geometric,arrows.meta,calc}
\usetikzlibrary{arrows, decorations.markings}
\definecolor{mydeepgreen}{HTML}{017A79}
\usetikzlibrary{patterns}
\allowdisplaybreaks
       \newtheorem{theorem}{\bf Theorem}[section]
         \newtheorem{conjecture}{\bf Conjecture}[section]

       \newtheorem{definition}{\bf Definition}[section]
       \newtheorem{remark}{\bf Remark}[section]

       \numberwithin{equation}{section}

\definecolor{mydeepgreen}{RGB}{3,100,50}

\begin{document}
\title{{\sl Weighted Birkhoff averages: Deterministic and probabilistic perspectives}}
\author{Zhicheng Tong $^{\mathcal{z}}$, Yong Li $^{\mathcal{x}}$}

\renewcommand{\thefootnote}{}
\footnotetext{\hspace*{-6mm}
\begin{tabular}{l   l}	$^\mathcal{z}$~School of Mathematics, Jilin University, Changchun 130012, P. R.  China. \url{tongzc25@jlu.edu.cn}\\	$^{\mathcal{x}}$~The corresponding author. School of Mathematics, Jilin University, Changchun 130012, P. R.  China; \\  Center for Mathematics and Interdisciplinary Sciences,  Northeast  Normal  University, Changchun 130024, \\P. R. China.  \url{liyong@jlu.edu.cn}
\end{tabular}}

\date{}
\maketitle

\begin{abstract}
In this paper, we survey physically related applications of a class of weighted quasi-Monte Carlo methods from a theoretical, deterministic perspective, and establish quantitative universal rapid convergence results via various regularity assumptions. Specifically, we introduce weighting with compact support to the Birkhoff ergodic averages of quasi-periodic, almost periodic, and periodic systems, thereby achieving universal rapid convergence, including both arbitrary polynomial and exponential types. This is in stark contrast to the typically slow convergence in classical ergodic theory. As new contributions, we not only discuss more general weighting functions but also provide quantitative improvements to existing results; the explicit regularity settings facilitate the application of these methods to specific problems. We also revisit the physically related problems and, for the first time, establish universal exponential convergence results for the weighted computation of Fourier coefficients, in both finite-dimensional and infinite-dimensional cases. In addition to the above, we explore results from a probabilistic perspective, including the weighted strong law of large numbers and the weighted central limit theorem, by building upon the historical results. 
\end{abstract}

\tableofcontents

\section{Introduction}
\renewcommand{\thefootnote}{\fnsymbol{footnote}}
For any given $ p,q>0 $, define a special normalized weighting function with compact support as follows:
\begin{equation}\label{CLTFUN}
	{w}_{p,q}\left( x \right):= \left\{ \begin{array}{ll}
		{\left( {\int_0^1 {\exp \left( { - {s^{ - p}}{{\left( {1 - s} \right)}^{ - q}}} \right)ds} } \right)^{ - 1}}\exp \left( { - {x^{ - p}}{{\left( {1 - x} \right)}^{ - q}}} \right), & \quad x \in \left( {0,1} \right),  \hfill \\
		0,&\quad x \notin \left(0,1\right). \hfill \\
	\end{array}  \right.
\end{equation}
It is observed that  $ {w_{p,q}}\left( x \right) \in C_0^\infty \left( {[0,1]} \right) $. Under \eqref{CLTFUN},  the corresponding weighted Birkhoff average  evaluated along a quasi-periodic/almost periodic trajectory of length $ N $ is given by (the discrete mapping case or the continuous flow case, respectively)
\begin{equation}\label{WBA}
	{\rm{WB}}_{N}\left( f \right)\left( \theta \right):=\frac{1}{{{A^{p,q}_N}}}\sum\limits_{n = 0}^{N - 1} {w_{p,q}\left( {n/N} \right)f\left( {{\mathscr{T}^n}\theta} \right)}\quad \text{or} \quad \frac{1}{N}\int_0^N {w_{{p,q }}\left( {s/N} \right)f\left( {{\mathscr{T}^s}\theta } \right)ds} ,
\end{equation}
where $ {A^{p,q}_N}: = \sum\nolimits_{n = 0}^{N - 1} {w_{p,q}\left( {n/N} \right)}  $ is a standard normalizer, $ f $ is a sufficiently smooth observable, $ \mathscr{T} $ is an irrational translation on the $ d $-torus $ \mathbb{T}^d :={\mathbb{R}^d}/{ {\mathbb{Z}} ^d}$ with $ d \in \mathbb{N}^+ \cup \{  + \infty \}  $\footnote{Here, $ \mathbb{T}^\infty $ could represent either the   infinite-dimensional torus  $ \mathbb{T}^{\mathbb{N}} $ or the   bi-infinite-dimensional torus  $ \mathbb{T}^\mathbb{Z} $. The same applies to $ \mathbb{R}^\infty $ and $ \mathbb{Z}^\infty $ in the subsequent text.}, and $ \theta\in \mathbb{T}^d $ is an initial value. More precisely, $ \mathscr{T}^s\theta= \theta +s\rho \mod 1 $ in each coordinate with some nonresonant (rationally-independent) rotation $ \rho \in \mathbb{R}^d $. We will illustrate this weighting procedure later in Figure \ref{fig:main}. It is worth emphasizing that, despite the fact that such quasi-periodic and almost periodic systems are classical and relatively simple, many fundamental nonlinear systems in dynamical systems and PDEs are always smoothly conjugated to them (e.g., via the Kolmogorov-Arnold-Moser (KAM) theory, etc). Therefore, it is \textit{natural} that, in the literature, researchers have chosen to study them directly in weighted Birkhoff averages of the form \eqref{WBA}.

This weighted approach, initially introduced as a breakthrough, notably in \cite{Las93a,Las93b,Las99}, exhibits unexpectedly excellent convergence, known as ``super convergence'', in the study of celestial mechanics as well as dynamical systems. This has been observed in practical computations \cite{Las93a,Las93b,Las99,DSSY17,DY18,SY18,SM20,MS21,BC23,DM23,BM24,Boh24,CCGd24,BdJC25,BHS25,BBC25,GZE25,MS25,SM25b,RB25,RKB25,SM25a,TL25b,ZXH25} and rigorously proved in theory \cite{DSSY17,DY18,DM23,RB24,TL24a,TL24b,TL25a,TL25c} as \textit{arbitrary polynomial} uniform convergence in the finite-dimensional setting. Here, we would like to mention that the original proof came from \cite{DSSY17,DY18}\footnote{These works originate from the research group of Prof. James Alan Yorke.}. In particular, the uniform convergence could even be \textit{exponential} \cite{TL24a,TL24b,TL25a,TL25b,TL25c} in a \textit{universal} sense, even in the \textit{infinite-dimensional} setting \cite{TL24b}. More precisely, the uniformity is for all initial values, while the universality is for almost all rotations with  full Lebesgue measure or full probability measure, as well as  observables that possess analyticity. We will provide later a detailed, rigorous and \textit{improved} statement in Theorem \ref{QMC}. However, it is worth mentioning that under this weighted method, the $ C^\infty_0([0,1]) $ property of the weighting function is \textit{essential}. Otherwise, for $ C^m_0([0,1]) $ weighting functions with finite $ m $, only polynomial convergence of finite order can be achieved in general, as shown by the counterexamples constructed in \cite{TL24b,TL25a,TL25c}.  This means that the weighting functions initially and frequently utilized, such as $  \sin^2(\pi x) $ \cite{Las93a,Las93b,Las99}, offer no advantages compared to \eqref{CLTFUN} in practice\footnote{From the perspective of ergodic theory, which does not focus solely on the rate of numerical simulation, considering the corresponding convergence results remains independently interesting.}, especially since they have essentially the same computational cost as \eqref{CLTFUN}, as pointed out in \cite{DSSY17,SM25a} among others.

However, in the classical ergodic theory, Birkhoff averages, with uniformly equal weights $ N^{-1} $, 		 depressingly possess \textit{quite slow} convergence.  Given a mapping $\mathcal{T}$ on a topological space $ X $ with a  probability measure $ \mu $ for which $ \mathcal{T} $ is invariant and ergodic, the Birkhoff average of an observable $ f $ defined on $ X $ is defined as
\[{{\rm B}_N}\left( f \right)\left( x \right): = \frac{1}{N}\sum\limits_{n = 0}^{N - 1} {f\left( {{\mathcal{T}^n}x} \right)} , \quad x \in X.\]
As the time length $N$ increases,  the Birkhoff ergodic theorem states that $ {{\rm B}_N}\left( f \right)\left( x \right) $ converges to the spatial average $ \int_X {fd\mu }  $ in a suitable way. Specifically, this convergence holds in the $L^2$ sense when $f$ admits $L^2$ integrability (von Neumann's ergodic theorem), and in the almost everywhere sense when $f$ admits $L^1$ integrability (Birkhoff's ergodic theorem).
As is well known since  \cite{Kac96}, for non-trivial (non-identically constant) observables, the convergence rate of the Birkhoff average is \textit{at most} $ \mathcal{O}(N^{-1}) $, and this situation typically occurs only in special cases.
It was initially shown in \cite{Kre78} that one can construct a non-trivial continuous observable under which the Birkhoff average, both in the almost everywhere sense and the norm sense, admitted an \textit{arbitrarily slow} convergence. More precisely, given any ergodic invertible measure preserving transformation $ \mathcal{T} $ of the interval $ [0,1] $ and any slowly varying sequence $0< \epsilon_N \to0^+ $, there exists a continuous observable $ f $ 
such that
\[\mathop {\lim \sup }\limits_{N \to  + \infty } \epsilon _N^{ - 1}\left| {{{\rm B}_N}\left( f \right)\left( x \right) - \int_X {fd\mu } } \right| =  + \infty , \quad \text{a.e.}\]
and
\[\mathop {\lim \sup }\limits_{N \to  + \infty } \epsilon _N^{ - 1}{\left\| {{{\rm B}_N}\left( f \right)\left( x \right) - \int_X {fd\mu } } \right\|_{{L^p (X)}}} =  + \infty. \]
This important observation is commonly referred to as ``\textit{the absence of estimates for rates of convergence in the Birkhoff ergodic theorem}''.
Very recently, similar yet distinct statements were proved in \cite{Ryz23,Ryz25a,Ryz25b,Ryz25c}\footnote{Z. Tong expresses deep gratitude for the communications with Prof. Valery V. Ryzhikov. He pointed out that similar results can be achieved by introducing weighting to the Birkhoff averages in the general case.} via elegant  approaches.             We would also like to mention other fundamental results on Birkhoff averages and related ergodic problems with slow convergence rates, for instance, \cite{dR79,KP81,PUZ89,LOT99,FS03,KT03,KP16,FS18,Fan19,KP19,HSY19a,HSY19b,Fan21,HLW21,Pod22,Pod24a,Pod24b,Ryz24} and the literature cited therein.

The general situation is not optimistic, even the classical case  of finite-dimensional toral translation  is no exception. These points can be seen more clearly here. Consider the previous Birkhoff average $ {{\rm B}_N}\left( f \right)\left( \theta \right) $, with $\theta\in  X=\mathbb{T}^d$ ($ d \in \mathbb{N}^+ $), $ \mathcal{T}=\mathscr{T}: \theta \to \theta+\rho $, and $ \mu $ being the Haar measure induced by the Lebesgue measure. On the one hand, as previously mentioned, $ {{\rm B}_N}\left( f \right)\left( \theta \right) $ can achieve its maximum possible rate of $\mathcal{O}^{\#}(N^{-1})$ only under special circumstances. This is one such case:  for almost all rotations $\rho \in \mathbb{R}^d$, if the observable $f$ is sufficiently smooth. This is evident by utilizing the co-homological equation argument. More precisely, for the given observable $ f :\mathbb{T}^d \to \mathbb{R}^1$ with zero average, let us consider the so-called co-homological equation
\begin{equation}\label{TDFC}
	f\left( \theta \right) = g\left( {\mathscr{T}\theta} \right) - g\left( \theta \right), \quad \theta \in {\mathbb{T}^d}.
\end{equation}
It admits a unique (up to a constant) solution $ g:\mathbb{T}^d \to \mathbb{R}^1 $, often referred to as the co-boundary,  with its Fourier expansion explicitly given by
\[g\left( \theta \right) = \sum\limits_{ 0 \ne  k \in {\mathbb{Z}^d}} {\frac{{\widehat f\left( k \right)}}{{\exp \left( {2\pi i\left\langle {k,\rho } \right\rangle } \right) - 1}} \cdot \exp \left( {2\pi i\left\langle {k,\theta} \right\rangle } \right)} , \quad \theta \in {\mathbb{T}^d}.\]
With the universal Diophantine nonresonance for $ \rho\in \mathbb{R}^d $ (see Definition \ref{Finite-dimensional Diophantine nonresonance}), it can be proved that $ g $ is well-defined on $ \mathbb{T}^d $, whenever $ f \in C^{L}(\mathbb{T}^d) $ with $ L>0 $ sufficiently large. For more accurate statements, see \cite{Kat04}. Therefore, by summing and averaging the co-homological equation \eqref{TDFC}, the convergence rate of $ {{\rm B}_N}\left( f \right)\left( \theta \right) $ can be obtained explicitly as $\mathcal{O}^{\#}(N^{-1})$.		
On the other hand, results with slow convergence rates can also be explicitly constructed. To show that there is no analogue of the so-called Denjoy-Koksma inequality in the case of $ \mathbb{T}^2 $, counterexamples were constructed in \cite{Yoc80,Yoc95a} using extremely Liouville rotations (see Section \ref{SUBSEC111} for the definition) and obtained \textit{arbitrarily slow} convergence via analytic observables. Unlike the case in \eqref{WBA}, for uniform weights---that is, weights that are independent of the time scale $ N $---one often can only obtain a convergence rate slower than $ \mathcal{O}^{\#}(N^{-1}) $, even when considering the case of toral translation  \cite{FS18,Fan19}.

Clearly, the slow convergence results mentioned above are fundamentally different from the universal arbitrary polynomial and exponential convergence results exhibited by \eqref{WBA} under non-uniform weighting. As one of the fundamental cornerstones of dynamical systems, ergodic theory, faces the challenging question of investigating its convergence rates in various settings. Although significant research has been conducted on weighted Birkhoff averages \eqref{WBA}, many essential issues regarding convergence rates remain unexplored. Our goal is not only to provide an almost complete theoretical guarantee for the practical computation of \eqref{WBA}, but more importantly, \textit{to explore its intrinsic connections with classical ergodic theory}. From a weighted perspective, many new and surprising phenomena arise, which are closely related to ergodic theory, number theory, and harmonic analysis. Consequently, this has independent interest.

Therefore, the \textbf{key motivations} of this paper are elaborated below.
\begin{itemize}
\item [(M1)] Given that the current results on arbitrary polynomial and exponential convergence, which are most suitable for application in the full Lebesgue measure or full probability measure sense and with analytic observables, are mostly presented in the form of corollaries and are highly scattered, we aim to \textit{integrate} these results and \textit{extend} them to arbitrary $ p,q>0 $  to facilitate their application in celestial mechanics and dynamical systems, see Theorem \ref{QMC}.

\item [(M2)] The vast majority of existing results focus on the special case in \eqref{WBA} where $p=q=1$, and are more inclined towards qualitative exponential convergence. However, there is currently \textit{no} quantitative universal exponential convergence for the infinite-dimensional case---which is more \textit{challenging} than the finite-dimensional case---regardless of any values of $p$ and $q$. We will solve this important problem in Theorem \ref{QMC}.

\item [(M3)] There have been numerous studies on quasi-periodic and almost periodic cases for \eqref{WBA}. However, the universal exponential convergence in the periodic case---which intuitively should be more straightforward---has \textit{yet} to be established. We will supplement this result \textit{quantitatively} in Theorem \ref{QMC}.

\item [(M4)] Regarding the weighted convergence rates in \eqref{WBA} with other similar weighting functions, there has been \textit{significant interest} in how fast universal exponential convergence can be achieved, and what the specific form might be. Theorems \ref{QMC} and \ref{TH3}, together with Remark \ref{RE13}, will discuss this issue.

\item [(M5)] One important application of \eqref{WBA} is the efficient computation of Fourier coefficients. However, there is currently \textit{no} rigorous proof of its exponential convergence, and discussions regarding the almost periodic case are also \textit{lacking}. In addition to these aspects, the \textit{effective order} should also be considered. It indicates the range within which the Fourier coefficients can be calculated at an exponentially fast rate. We will give an almost complete resolution of this problem in Theorem \ref{CLT-FOU}.

\item [(M6)] There have already been many results in the deterministic case under the weighting function \eqref{CLTFUN}, and we also hope to further explore results from a \textit{probabilistic} perspective as well. The primary consideration behind this motivation is that this weighted form bears similarities to the Gaussian distribution in statistics.  We will present several  results in Theorems \ref{SJT1} to \ref{CLTT1}.
\end{itemize}

Based on the aforementioned various motivations, the \textbf{main contributions} and organization of this paper are as follows. In Section \ref{SEC11}, we review, summarize and \textit{extend} the results of the weighted method in \eqref{WBA}---alternatively referred to as a \textit{weighted quasi-Monte Carlo method}---from a  deterministic perspective. In Section \ref{SUBSEC111}, we first introduce a limit-preserving  theorem and its converse (Theorem \ref{LIMIT}), which is essentially an abstract averaging version of the Toeplitz theorem. This theorem allows us to directly employ such weighting to accelerate convergence, provided that the limit of the average exists \textit{in advance}. In addition to this, we also introduce some fundamental concepts in both finite-dimensional and infinite-dimensional settings, such as the spatial structure, Diophantine nonresonance, and analyticity. Next, as one of the \textit{main} results of this paper, Theorem \ref{QMC} provides convergence rate results for \eqref{WBA} under various regularity conditions---especially analyticity---in a universal sense with full Lebesgue measure or full probability measure. These results cover quasi-periodic, almost periodic, and periodic cases. The majority of the results are presented in the form of \textit{quantitative} exponential convergence and are applicable for any $p,q > 0$. In particular, the results for the almost periodic and periodic cases are \textit{entirely new}. These findings not only provide \textit{nearly complete} theoretical guarantees for weighted fast computation in celestial mechanics and dynamical systems but also highlight the connections with ergodic theory, number theory, and harmonic analysis. As a new discovery, we introduce a novel weighting function in Theorem \ref{TH3}, which enables the weighted Birkhoff average to converge at a rate that is \textit{nearly linear-exponential}. All the above results are stated in the toral translation form for brevity, as explained previously. They could also be restated in the \textit{conjugated form}  as those in \cite{DSSY17,DY18,DM23,RB24}, which we do not pursue here. In Section \ref{SECphy}, we revisit at least five applications of this weighted method in physically related problems. As one of the important applications, it can compute Fourier coefficients with an \textit{exponential} convergence rate. In Section \ref{SECFOU}, we rigorously derive the corresponding approximation formula in a weighted manner and, \textit{for the first time}, provide a proof of the exponential convergence of the weighted average and determine the \textit{effective order}, including both finite-dimensional and infinite-dimensional cases. This result, which is \textit{not} a direct consequence of Theorem \ref{QMC}, significantly differs from existing results and addresses substantial challenges, thus constituting a \textit{primary} contribution of this paper. We also briefly discuss, in Section \ref{SEC114}, which  factors may affect the ``practical'' convergence rate, thereby providing a more dynamical-systems perspective for numerical simulations. Apart from the deterministic results, we also aim to explore the probabilistic perspectives in Section \ref{SECPRO}. Since the various weighted laws of large numbers are almost trivial at this point, we choose not to state all of them specifically. Based on historical results, Section \ref{SECWL} presents the weighted strong laws of large numbers (Theorems \ref{LIMT1} and \ref{LIMT2}); Section \ref{SECWCLT} establishes the weighted central limit theorem (Theorem \ref{CLTT1}) and leaves a conjecture open. To clarify the new contributions of this paper and to compare them with the historical literature, we add Section \ref{SECCON} as a concise summary for  readers. Finally, Section \ref{SECPD} provides the proof of the deterministic results, while Section \ref{PPR} contains the proof of the probabilistic results.

\subsection{The deterministic perspective}\label{SEC11}
\setcounter{footnote}{0}
\renewcommand{\thefootnote}{\fnsymbol{footnote}}

\subsubsection{The weighted quasi-Monte Carlo method}\label{SUBSEC111}
As one of the most common numerical integration methods, the renowned Monte Carlo method estimates the value of integrals through random sampling; it indeed has extensive applications in both deterministic and probabilistic problems  \cite{Sha10}. It is  well-suited for calculating integrals that are difficult to handle using traditional analytical methods, especially those involving high dimensions and complex regions. However, its application involves sampling from intricate, high-dimensional probability distributions, a process that is not only challenging but also demands significant computational resources. In addition, the convergence rate of the Monte Carlo method is generally slow. Therefore, researchers have sought alternative approaches to improve it, collectively known as \textit{quasi-Monte Carlo methods}. 

\begin{figure}[htbp]
 	\centering
 	 
 	\begin{subfigure}[b]{0.47\textwidth}
 		\centering
 		\includegraphics[width=\textwidth]{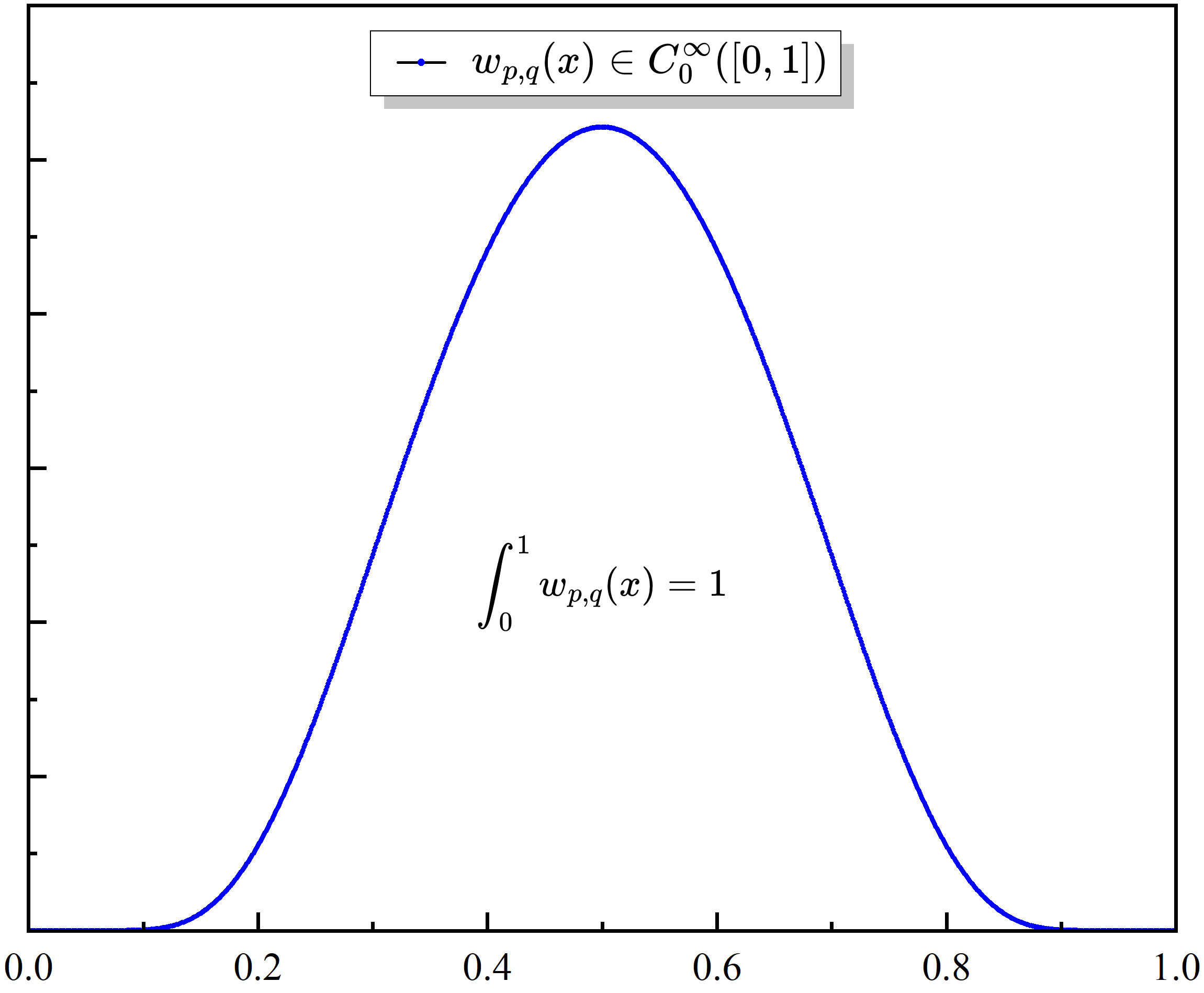}
 		\caption{The shape of the weighting function $ w_{p,q}(x) $ and its normalization property.}
 		\label{fig:sub1}
 	\end{subfigure}
 	\hfill 
 	\begin{subfigure}[b]{0.47\textwidth}
 		\centering
 		\includegraphics[width=\textwidth]{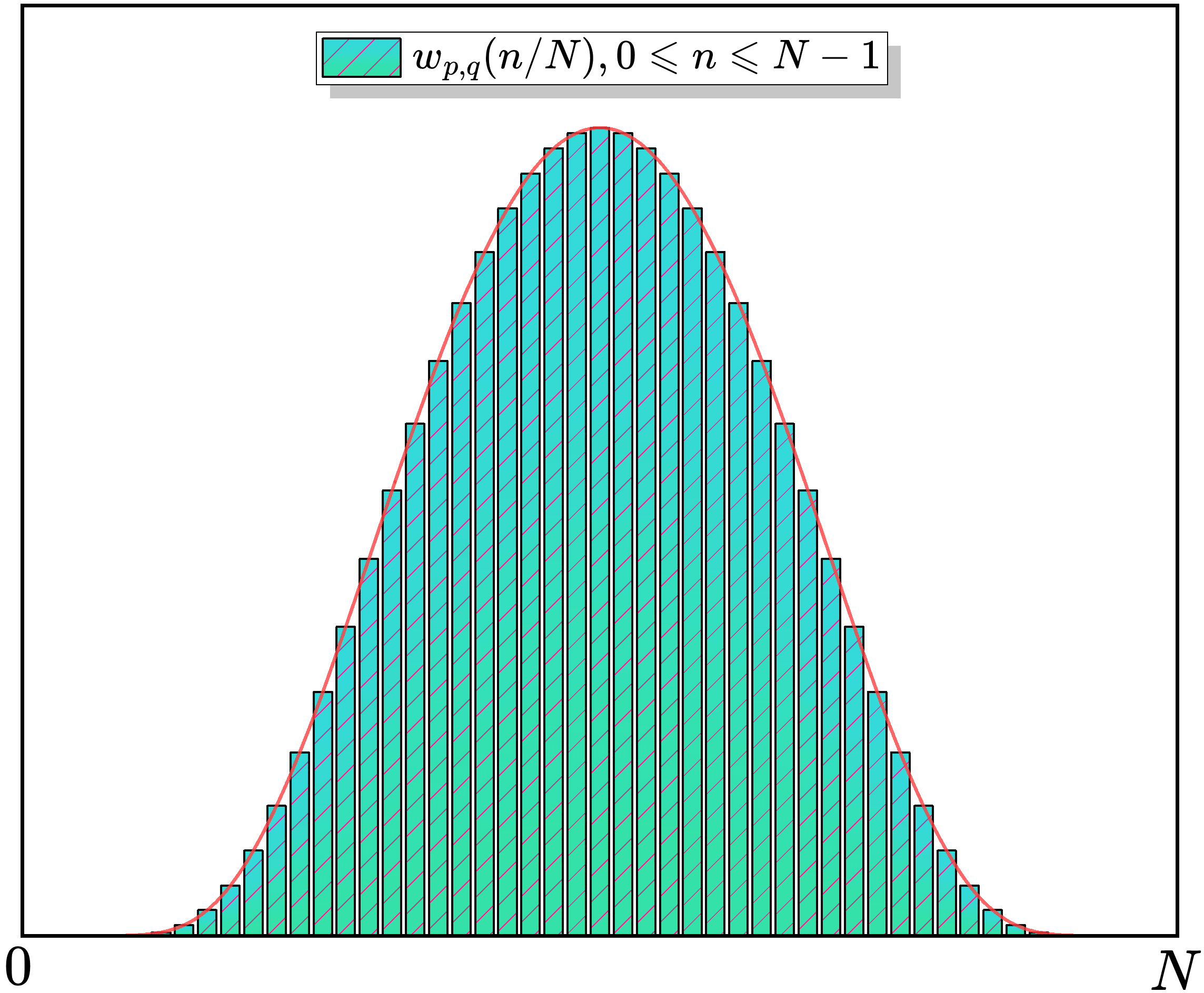}
 		\caption{Uniformly sampled weights from the weighting function $ w_{p,q}(x) $ over the time scale $ [0, N-1] $.}
 		\label{fig:sub2}
 	\end{subfigure}
 	\par\bigskip  
 	\begin{subfigure}[b]{0.47\textwidth}
 		\centering
 		\includegraphics[width=\textwidth]{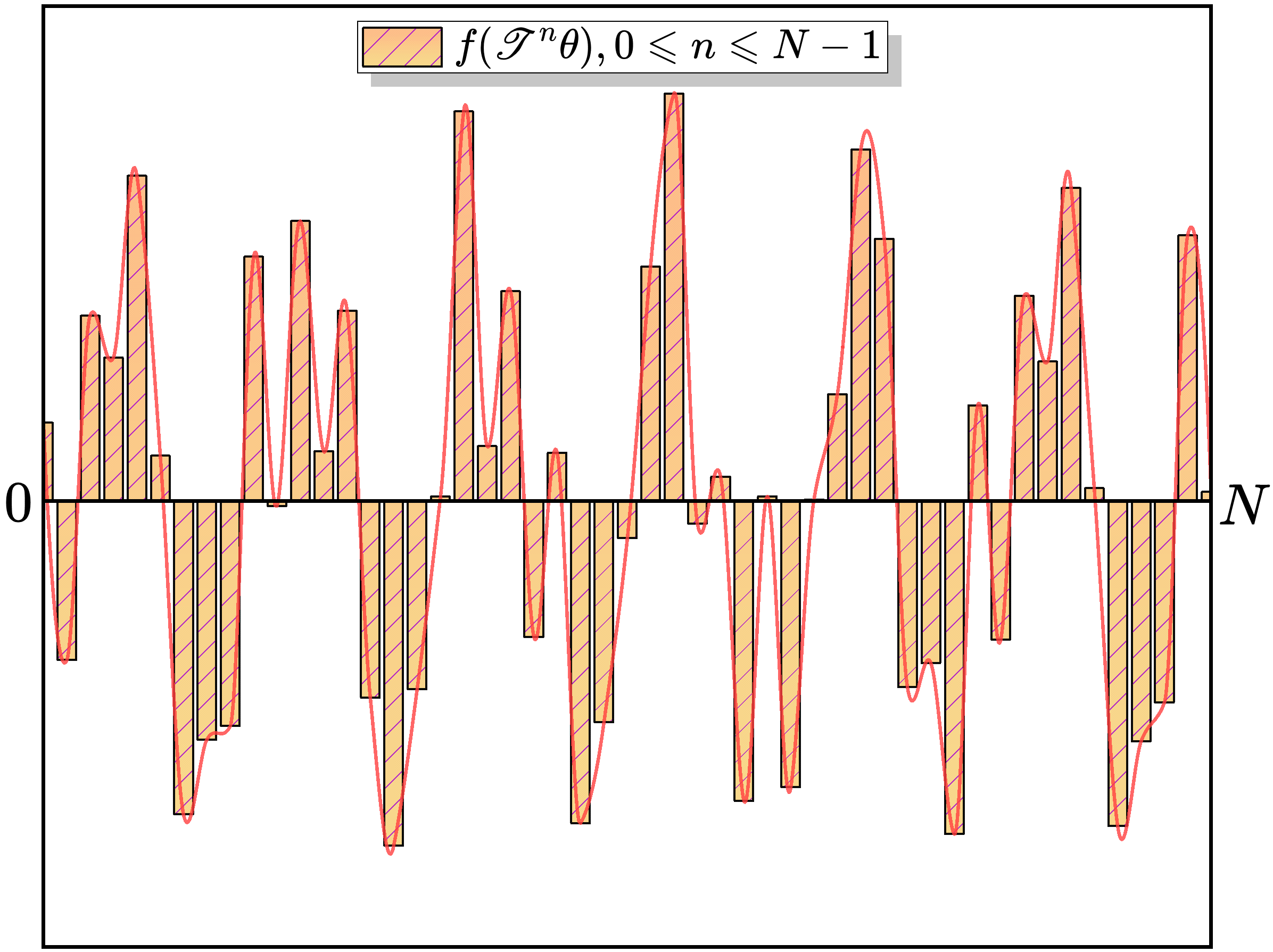}
 		\caption{Some data generated based on a quasi-periodic dynamical system, with the mean assumed to be zero here.}
 		\label{fig:sub3}
 	\end{subfigure}
 	\hfill  
 	\begin{subfigure}[b]{0.47\textwidth}
 		\centering
 		\includegraphics[width=\textwidth]{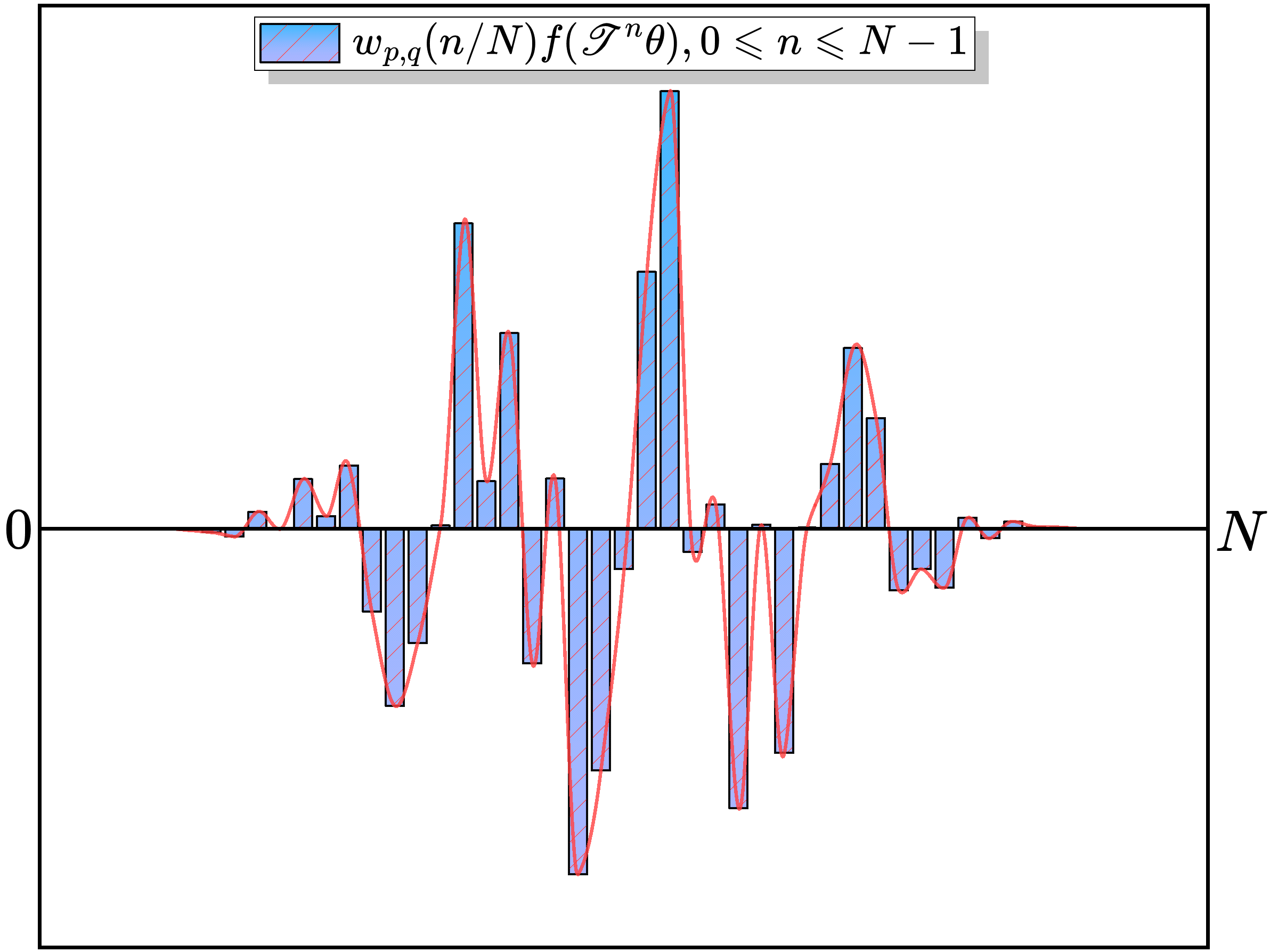}
 		\caption{The data in \eqref{fig:sub3} are weighted using the weights depicted in \eqref{fig:sub2}, which significantly reduces the weights of the initial and terminal data points.}
 		\label{fig:sub4}
 	\end{subfigure}
 	\caption{The weighting procedure of the weighted quasi-Monte Carlo method for data generated from a quasi-periodic dynamical system}
 	\label{fig:main}
 \end{figure}

The weighted method presented in this paper, as illustrated in \eqref{WBA}, is also a quasi-Monte Carlo technique and can be appropriately termed the \textit{(non-uniformly) weighted quasi-Monte Carlo method}. Here, the term ``non-uniformly'' means that the weights depend on the time scale, and we would like to omit it in the following discussions for brevity.
 As previously mentioned, it significantly accelerates calculations involving averages (or integrals) in many  dynamical systems of interest.   Intuitively speaking, based on the weighted approach and the excellent compact support property of the weighting function $  w_{p,q}(x) $, this weighted quasi-Monte Carlo method significantly reduces the weight of the initial and terminal data, focusing more on the intermediate data. The elaboration of the weighting procedure is shown explicitly in Figure \ref{fig:main}. However, it is important to note that,  the intermediate data is not necessarily crucial either, and the weight of a significant portion of the data can also be substantially reduced (i.e., by utilizing some new weighting functions) while still ensuring rapid convergence of arbitrary polynomial or even exponential types \cite{TL25a}. Therefore, we can actually allow the data to be \textit{sparse}.

 When applying the weighting function $ w_{p,q} (x)$ to the original average, it is necessary to ensure that it does not change the original limit. Without this assurance, this weighted method may not be significant, particularly in terms of accelerating the convergence rate. The theorem below essentially represents an abstract \textit{averaging version of  the Toeplitz theorem}. Dealing with the Banach space $ \mathbb{R} $ endowed with the sup-norm, it was claimed in \cite{DSSY17} this assertion   without a proof. In this paper, we provide  a rigorous and concise  proof for this intuitive phenomenon, as detailed  in Section \ref{PROOFPRO1}.

 \begin{theorem}\label{LIMIT}
 	Consider  a sequence $ \{a_n\}_{n=1}^{\infty}$ that  belongs to some Banach space $ (\mathscr{B},\|\cdot\|_\mathscr{B}) $. Assume that 
 	\begin{equation}\label{OPTpro1}
 		\frac{1}{N}\sum\limits_{n = 1}^N {{a_n}}  \to a \in \mathscr{B} \quad \text{in the $ \|\cdot\|_\mathscr{B}$ norm.}
 	\end{equation}
 	Then we have
 	\begin{equation}\label{OPTPro2.1}
 		\frac{1}{{{A^{p,q}_N}}}\sum\limits_{n = 1}^{N} {w_{p,q}\left( {n/N} \right){a_n}}  \to a \in \mathscr{B} \quad \text{in the $ \|\cdot\|_\mathscr{B}$ norm.}
 	\end{equation}
  The same conclusion holds for the  continuous case (the integral version). However, the converse is not true, i.e., \eqref{OPTPro2.1} does not imply \eqref{OPTpro1}. In fact, similar statements hold true for all weighting functions in $ C_0^1([0,1]) $.
 \end{theorem} 
 \begin{remark}
	Even if the convergence of the original average \eqref{OPTpro1} is assumed, this weighted quasi-Monte Carlo method may not always significantly increase the convergence rate. For instance, if $ a_n=n^{-1} $, it is evident that \eqref{OPTpro1} converges at a rate of $ \mathcal{O}^{\#}(N^{-1} \log N  ) $, while \eqref{OPTPro2.1} converges at a rate of  $ \mathcal{O}^{\#}(N^{-1}) $, which only eliminates the logarithmic factor. See also Section \ref{SECphy} for more cases.
\end{remark}

The primary significance of Theorem \ref{LIMIT} is to demonstrate that anything observed using the law of large numbers can also be observed using this weighted quasi-Monte Carlo method. Theorem \ref{LIMIT} also tells us that even if the weighted average \eqref{OPTPro2.1} converges, the original one \eqref{OPTpro1} does not necessarily converge. Therefore, in numerical simulations \cite{Las93a, Las93b, Las99, DSSY17, DY18, SY18,   SM20, MS21, BC23, DM23, BM24,Boh24,CCGd24,BdJC25,BHS25, BBC25, GZE25,   MS25, RB25, RKB25, SM25a,  SM25b,TL25b,ZXH25}, the theoretical convergence of the average \eqref{OPTpro1} must be established \textit{in advance} (e.g., via ergodicity) before considering weighted acceleration.

 We next present a theoretical result (Theorem \ref{QMC}) that demonstrates the universal arbitrary polynomial and exponential convergence for quasi-periodic, almost periodic, and periodic dynamical systems via this weighted quasi-Monte Carlo method, for \textit{all} $ p,q>0 $. 
 The results for quasi-periodic and almost periodic systems have been obtained under various settings in \cite[Theorem 3.1]{DSSY17}, \cite[Theorem 1.1]{DY18},  \cite[Theorem 3]{DM23}, \cite[Theorem 2.1, Corollaries 2.1, 3.1 and 3.2]{TL24a}, \cite[Theorem 4.4]{TL24b}, and \cite[Theorems 6.1 and 6.2]{TL25b}, mainly for $ p=q=1 $\footnote{Here, we emphasize that the analysis for exponential convergence becomes more challenging when $ p,q>0 $ in general, as opposed to the simpler case where $ p = q = 1 $.}. Specifically, the result concerning the universal exponential convergence in the almost periodic case is  \textit{significantly improved} here, with the exponent explicitly depending on the spatial structure. The exponential convergence for the periodic case is presented here for the first time, although it can be regarded as a variation of \cite[Theorem 3.5]{TL24b} (if we ignore the quantitative aspects), and we would also like to mention the very recent work \cite[Theorem 2.2]{RB24} on arbitrary polynomial convergence for the periodic case with $ C_0^\infty ([0,1]) $ weighting functions. For the sake of simplicity, we prefer not to present such results with extensions to $ p,q>0 $ in general Banach spaces here. Interested readers are referred to \cite{TL24a,  TL24b,TL25a} for further details.
 
Before stating, we need to introduce some basic notions. Consider a vector-valued (could be infinite-dimensional, here and henceforth) observable $ f=(f_1,\ldots,f_m) $ on $ \mathbb{T}^d $ with $   d,m \in \mathbb{N}^+ \cup \{  + \infty \} $, as usual, its  $ \ell^\infty $-norm is then  given by $ \|  f \|_{\ell^\infty}:=  \mathop {\sup }\nolimits_m \left| {{f_m}} \right|$. If the observable is sufficiently smooth, then we have the corresponding Fourier analysis. It should be emphasized that the finite-dimensional and infinite-dimensional cases are quite different, hence they need to be considered separately.

	For $ d \in \mathbb{N}^+ $ and $ \sigma>0 $, the finite-dimensional thickened  torus $ \mathbb{T}_\sigma ^d $ is usually defined as
\begin{equation}\label{FDTT}
	\mathbb{T}_\sigma ^d : = \left\{ {x = {{({x_j})}_{1 \leqslant j \leqslant d}}: \quad {x_j} \in \mathbb{C},\;\operatorname{Re} {x_j} \in \mathbb{T}^1,\;\left| {\operatorname{Im} {x_j}} \right| \leqslant \sigma },\;\text{for all}\; 1 \leqslant j \leqslant d \right\}.
\end{equation}
By definition, $ \mathbb{T}_\sigma ^d $ amounts to a complex analytic continuation of the standard torus $ \mathbb{T} ^d={\mathbb{R}^d}/{\mathbb{Z}^d} $, see \cite{Sal04} for details. To provide readers an intuitive grasp of the $ d=1 $ case, we present a schematic illustration in Figure \ref{fig:torus-strip} below.
\begin{figure}[htbp]
	\centering \begin{tikzpicture}[scale=2.0]
	\draw[line width=0.8pt,->] (-1.2,0)--(2.2,0) node[below] {$x$};
	\draw[line width=0.8pt,->] (0,-1.2)--(0,1.2) node[left]  {$y$};
	
	\def\mysigma{0.6}
	
\draw[line width=1.5pt,mydeepgreen,
pattern=north east lines, pattern color=mydeepgreen!60]
(0,-\mysigma) rectangle (1,\mysigma);
	
	\draw[line width=1pt,mydeepgreen] (0,-\mysigma) rectangle (1,\mysigma);

	\draw[line width=1.5pt, red] (0,0)--(1,0);
	
	\node[below=3pt,red] at (0.5,0) {$\mathbb{T}^1$};
	
		\node[below=3pt,mydeepgreen] at (1.4,0.6) {$\mathbb{T}_\sigma^1$};
	
		\node[below=3pt,black] at (1.8,1.3) {$\mathbb{C}^{1}:z=x+iy$};
		
		\fill (1,0) circle (1pt);
		\node[below=3pt] at (1.4,0) {$x=1$};

		\fill (0,\mysigma) circle (1pt);
		\node[left=3pt]  at (0,\mysigma) {$y=\sigma$};
		
			\fill (0,-\mysigma) circle (1pt);
		\node[left=3pt]  at (0,-\mysigma) {$y=-\sigma$};
		
			\fill (0,0) circle (1pt);
		
	\node[below left] at (0,0) {$O$};
\end{tikzpicture}
\caption{The shape of the  $ 1 $-dimensional thickened  torus $ \mathbb{T}_\sigma ^1 $ with $ \sigma>0 $}
\label{fig:torus-strip}
\end{figure}

With the definition of $\mathbb{T}_{\sigma}^d $, the finite-dimensional analyticity of observables is consequently given by:

\begin{definition}[Finite-dimensional analyticity]\label{Finite-dimensional analyticity}
	We say that the observable $ f $ is  analytic on $ \mathbb{T}^d $ with $ d \in \mathbb{N}^+ $, if the domain of $ f $ can be extended to $ \mathbb{T}_\sigma ^d $ with some $ \sigma>0 $\footnote{The largest $ \sigma>0 $ is usually called the analytic radius of $ f $. Similar expressions hold in the infinite-dimensional case as well.}, and 
	\[f  = \sum\limits_{k  \in \mathbb{Z}^d } {\widehat f (k) \exp\left({2\pi i \left\langle {k,x} \right\rangle  }\right)}\quad\text{with}\quad\|f\|_{\sigma,d }: = \sum\limits_{k  \in \mathbb{Z} ^d } {\| {\widehat f (k)} \|_{\ell^\infty}{\exp\left(2 \pi \sigma {{\| k  \|_{\ell^1}} }\right)}}  <  + \infty .\]
Here  $ \left\langle {k,x} \right\rangle : = \sum\nolimits_{j = 1}^d {{k_j}{x_j}}  $ denotes   the inner product, and $\left\|k\right\|_{\ell^1}:=\sum\nolimits_{j = 1}^d {\left| {{k_j}} \right|} $ denotes the $ \ell^1 $-norm  for all $ k \in \mathbb{Z}^d $.
\end{definition}

The finite-dimensional Diophantine nonresonance condition for rotations can be formulated as follows:

\begin{definition}[Finite-dimensional Diophantine nonresonance]\label{Finite-dimensional Diophantine nonresonance}
	We say that $ \rho \in \mathbb{R}^d $ is Diophantine, if there exists $ \gamma > 0 $ such that: 
	\begin{itemize}
		\item [(I)] The discrete case: for any $ 0\ne k \in \mathbb{Z}^d $ and any $ n \in \mathbb{Z} $, $ \left| {\left\langle {k,\rho } \right\rangle  - n} \right| > \gamma \left\| k \right\|_{{\ell ^1}}^{ - \tau } $, provided $ \tau\geqslant d $;
		
		\item [(II)] The continuous case: for any $ 0\ne k \in \mathbb{Z}^d $, $ \left| {\left\langle {k,\rho } \right\rangle } \right| > \gamma \left\| k \right\|_{{\ell ^1}}^{ - \tau } $, provided $ \tau\geqslant d-1 $.
	\end{itemize}
\end{definition}

Here, $ \tau $ is termed the Diophantine exponent of $ \rho $. We say that $ d $ in (I) and $ d-1 $ in (II) within Definition \ref{Finite-dimensional Diophantine nonresonance} are \textit{critical} Diophantine exponents. This is because the set of such rotations in $ \mathbb{R}^d $ admits zero Lebesgue measure but full Hausdorff dimension (e.g, $ \sqrt{2} $ and the golden ratio $ (\sqrt{5}+1)/2 $ for $ d=1 $). They do not exist when the Diophantine exponent is smaller, as can be proved by Dirichlet's principle. Conversely, when the Diophantine exponent is larger, the set of these rotations admits full Lebesgue measure and can thus be termed ``universal'' (e.g., $ \pi $ for $ d=1 $). For further details in this aspect, see \cite{Her79}. If a rotation is not Diophantine, we refer to it as Liouville.  Further classifications exist for Liouville rotations, such as  Brjuno (Bruno/Bryuno) rotations\footnote{We specifically consider Brjuno rotations that are strictly weaker than Diophantine in the sense of nonresonance (i.e., excluding the Diophantine case itself), thus belonging to the Liouville class.}, which  often essentially characterize the optimal arithmetic properties that are specific to various analytic dynamical systems in a sum/integral way \cite{Yoc95b,Yoc02}.

However, when considering the infinite-dimensional torus $ {\mathbb{T}^\infty } $ \cite{Bou05,Koz20,Koz21,MP21,BCGP24,CGP24}, it becomes necessary to impose some spatial structure---which may not be unique---to prevent the Fourier series expansions from blowing up, for instance. Here, we directly employ the concepts from \cite{Bou05}, which were systematically refined in \cite{MP21} and are primarily utilized to investigate quasi-periodic and almost periodic solutions in PDEs via KAM theory, including the nonlinear Schr\"odinger equation. This is because, once the corresponding toral translation conjugation theorem is established, the results of this paper can be directly applied to these PDEs.

For a fixed integer $ 2 \leqslant \eta \in \mathbb{N}^+ $ (here and henceforth), we define the following set of infinite integer vectors with finite support:
\begin{equation}\notag
	\mathbb{Z}_ * ^\infty : = \left\{ {k \in {\mathbb{Z}^\infty}: \quad {{\left| k \right|}_\eta }: = \sum\limits_{j \in \mathbb{N}} {{{\left\langle j \right\rangle }^\eta }\left| {{k_j}} \right|}  <  + \infty ,\;\left\langle j \right\rangle : = \max \left\{ {1,\left| j \right|} \right\}} \right\}.
\end{equation}
When $k \in  	\mathbb{Z}_ * ^\infty  $ is fixed, it is evident that  $ {k_j} \ne 0 $ only for finitely many indices $ j \in \mathbb{N} $. Such a metric as $ {{{\left| \cdot \right|}_\eta }} $ is quite essential in the infinite-dimensional case, as it determines certain boundedness of the summation in the analysis, both in KAM theory and in weighted ergodic theory discussed in this paper. Moreover, for $ \sigma>0 $, 
 the thickened infinite-dimensional torus $ \mathbb{T}_\sigma ^\infty $ is usually defined as 
\begin{equation}\notag 
	\mathbb{T}_\sigma ^\infty : = \left\{ {x = {{({x_j})}_{j \in \mathbb{N}}}: \quad {x_j} \in \mathbb{C},\;\operatorname{Re} {x_j} \in \mathbb{T}^1,\;\left| {\operatorname{Im} {x_j}} \right| \leqslant \sigma {{\left\langle j \right\rangle }^\eta } },\; \text{for all}\; j \in \mathbb{N} \right\}.
\end{equation}
 To help  readers appreciate how it differs from the finite-dimensional case (see for instance Figure \ref{fig:torus-strip}), Figure \ref{fig:torus-strip2} below displays the $ 1 $-dimensional torus corresponding to the $ j $-th component space of $ \mathbb{T}^\infty_\sigma $.
\begin{figure}[htbp]
	\centering \begin{tikzpicture}[scale=2.0]
		\draw[line width=0.8pt,->] (-1.2,0)--(2.2,0) node[below] {$x$};
		\draw[line width=0.8pt,->] (0,-1.2)--(0,1.2) node[left]  {$y$};
		
		\def\mysigma{0.6}
		
		\draw[line width=1.5pt,mydeepgreen,
		pattern=north east lines, pattern color=mydeepgreen!60]
		(0,-\mysigma) rectangle (1,\mysigma);
		
		\draw[line width=1pt,mydeepgreen] (0,-\mysigma) rectangle (1,\mysigma);

		\draw[line width=1.5pt, red] (0,0)--(1,0);
		
		\node[below=3pt,red] at (0.5,0) {$\mathbb{T}^1$};
		
		\node[below=3pt,mydeepgreen] at (1.4,0.6) {$\left(\mathbb{T}_\sigma^\infty\right)_j$};
		
		\node[below=3pt,black] at (1.8,1.3) {$\mathbb{C}^{1}:z=x+iy$};
		
		\fill (1,0) circle (1pt);
		\node[below=3pt] at (1.4,0) {$x=1$};

		\fill (0,\mysigma) circle (1pt);
		\node[left=3pt]  at (0,\mysigma) {$y=\sigma {\left\langle j \right\rangle }^\eta$};
		
		\fill (0,-\mysigma) circle (1pt);
		\node[left=3pt]  at (0,-\mysigma) {$y=-\sigma {\left\langle j \right\rangle }^\eta$};
		
		\fill (0,0) circle (1pt);
		
		\node[below left] at (0,0) {$O$};
	\end{tikzpicture}
	\caption{The shape of the $ j $-component space of the  infinite-dimensional thickened  torus $ \mathbb{T}_\sigma ^\infty $ with $ \sigma>0 $}
	\label{fig:torus-strip2}
\end{figure}

Given the aforementioned spatial structure, let us now introduce the infinite-dimensional analyticity as defined below.

\begin{definition}[Infinite-dimensional analyticity]\label{Infinite-dimensional analyticity} 	We say that the observable $ f $ is  analytic on $ \mathbb{T}^\infty $, if the domain of $ f $ can be extended to $ \mathbb{T}_\sigma ^\infty $ with some $ \sigma>0 $, and 
	\[f  = \sum\limits_{k  \in \mathbb{Z}_*^\infty } {\widehat f (k) \exp\left({2\pi i \left\langle {k,x} \right\rangle  }\right)}\quad\text{with}\quad\|f\|_{\sigma , \infty}: = \sum\limits_{k  \in \mathbb{Z}_* ^\infty } {\| {\widehat f (k)} \|_{\ell^\infty}{\exp\left(2 \pi \sigma {{| k  |_{\eta}} }\right)}}  <  + \infty .\]
\end{definition}
The infinite-dimensional Diophantine nonresonance condition for rotations can be stated as follows:

\begin{definition}[Infinite-dimensional Diophantine nonresonance]\label{Infinite-dimensional Diophantine nonresonance}
	We say that $ \rho \in \mathbb{R}^\infty $ is Diophantine, if there exists $ \gamma > 0 $ such that: 
	\begin{itemize}
		\item [(I)] The discrete case: for any $ 0\ne k \in \mathbb{Z}_*^\infty $ and any $ n \in \mathbb{Z} $, 
		$ \left| {\left\langle {k,\rho } \right\rangle  - n} \right| > \gamma {\left( {\prod\nolimits_{j \in {\mathbb{N} }} {(1 + {{\left\langle j \right\rangle }^\tau }{{\left| {{k_j}} \right|}^\tau })} } \right)^{ - 1}} $,
		provided $ \tau>1 $;
		
		\item [(II)] The continuous case: for any $ 0\ne k \in \mathbb{Z}_*^\infty $, 		$ \left| {\left\langle {k,\rho } \right\rangle  } \right| > \gamma {\left( {\prod\nolimits_{j \in {\mathbb{N} }} {(1 + {{\left\langle j \right\rangle }^\tau }{{\left| {{k_j}} \right|}^\tau })} } \right)^{ - 1}} $,
		provided $ \tau>1 $.
	\end{itemize}
\end{definition}

In contrast to the delicate results in finite dimensions, it seems that there is no known so-called ``0--1 law'' for nonresonance in infinite dimensions. However, the infinite-dimensional Diophantine nonresonance in Definition \ref{Infinite-dimensional Diophantine nonresonance} is universal. More precisely, under the standard probability measure $ \mathbb{P} $ induced by the product measure of the cube $ {\left[ {1,2} \right]^\mathbb{\infty}} $ (equivalent  to $ {\mathbb{T}^\infty } $), the set of these rotations admits full probability measure \cite{Bou05,MP21}. This is sufficient for us to establish the universal exponential convergence of weighted Birkhoff averages, via analytic observables. It should be emphasized that the universal nonresonance is not necessarily confined to the aforementioned Diophantine form. One can construct more general cases \cite{TL24b}, which we do not pursue here.

Now, we are in a position to state Theorem \ref{QMC}, for general $ p,q>0 $.
 \begin{theorem}\label{QMC}
 	Consider the weighted Birkhoff average in \eqref{WBA} utilizing the  weighted quasi-Monte Carlo method with $ p,q>0 $. It uniformly (with respect to the initial value $ \theta $) converges to the spatial average $ \int_{{\mathbb{T}^d}} {fdx}  $ in the $ \ell^\infty $-norm. Moreover, we have the following results:
 	\begin{itemize} 
 		\item [(I)] The quasi-periodic case:
 			\begin{itemize}
 			\item [(1)] If the observable $ f $ is $ C^\infty $ on $ \mathbb{T}^d $, then for almost all rotations $ \rho \in \mathbb{R}^d $, the convergence can be of  arbitrary polynomial type. More precisely, the error is $ \mathcal{O}\left(N^{-m}\right) $ for any $ m>0 $.
 			
 			\item [(2)]  If the observable $ f $ is analytic on $ \mathbb{T}^d $, then for almost all rotations $ \rho \in \mathbb{R}^d $, the convergence can be of   exponential type. More precisely, the  error is $ \mathcal{O}\left( {\exp ( { - {N^{\zeta_1 }}} )} \right) $ for any $ 0 < {\zeta _1} < {\left( {d + 1 + 1/\min \left\{ {p,q} \right\}} \right)^{ - 1}} $.

 				\item [(3)] If the observable $ f $ is ``sub-super analytic'' on $ \mathbb{T}^d $, namely $ \|\widehat f(k)\|_{\ell^\infty}=\mathcal{O}\left(\exp(-\|k\|_{\ell^1}^v)\right) $ for any $ v>1 $,  then for any  given $ 0<\kappa<1 $, when $ p,q $ are sufficiently large, for almost all rotations $ \rho \in \mathbb{R}^d $, the convergence  can be of   arbitrary exponential type. More precisely, the  error is $ \mathcal{O}\left( {\exp ( { - {N^{\kappa }}} )} \right) $.
 				
 				\item[(4)] If the observable $ f $ is a finite trigonometric polynomial on $ \mathbb{T}^d $, then for almost all rotations $ \rho \in \mathbb{R}^d $, the convergence can be of   exponential type. More precisely, the  error is $ \mathcal{O}\left( {\exp ( { - c{N^{\zeta_2 }}} )} \right) $ for some universal constant $ c>0 $, where $ {\zeta _2} = {\left( {1 + 1/\min \left\{ {p,q} \right\}} \right)^{ - 1}} $.
 			\end{itemize}
 		\item [(II)] The almost periodic case:
 			\begin{itemize}
 			\item [(5)] If the observable $ f $ is analytic on $ \mathbb{T}^\infty $, then for almost all rotations $ \rho \in \mathbb{R}^\infty $, the convergence can be of   exponential type. More precisely, the  error is $ \mathcal{O}\left( {\exp ( - {{(\log N)}^{{\zeta _3}}})} \right) $ for any $ 2 \leqslant {\zeta _3} < 1 + \eta $.  
 			
 			\item [(6)] If the observable $ f $ is ``super analytic'' on $ \mathbb{T}^\infty$, namely $ \|\widehat f(k)\|_{\ell^\infty}=\mathcal{O}\left(\exp(-\exp(|k|_\eta)\right) $,  then for almost all rotations $ \rho \in \mathbb{R}^\infty $, the convergence can be of   exponential type. More precisely, the error is $ \mathcal{O}\left( {\exp ( { - {N^{\zeta_4 }}} )} \right) $ for any $ 	0 < {\zeta _4} < {\left( { 1 + 1/\min \left\{ {p,q} \right\}} \right)^{ - 1}} $.

 		\end{itemize}
 	\item [(III)] The periodic case:	\begin{itemize}\item[(7)]
 			 Any $ T $-periodic ($ T \in \mathbb{N}^+ $) data converge exponentially to the periodic average through this  weighted quasi-Monte Carlo method for $ p,q>0 $. More precisely, for  $ f(n)=f(n+T) $ with $ n \in \mathbb{N} $ ($ f $ need not be continuous),
 		\[\left| {\frac{1}{{{A^{p,q}_N}}}\sum\limits_{n = 0}^{N - 1} {w_{p,q}\left( {n/N} \right)f\left( n \right)}  - \frac{1}{T}\sum\limits_{s = 0}^{T - 1} {f\left( s \right)} } \right|=\mathcal{O}\left( {\exp ( { - c{N^{\zeta_5 }}} )} \right)\]
 		holds for some universal constant $ c>0 $, where  $ {\zeta _5} = {\left( {1 + 1/\min \left\{ {p,q} \right\}} \right)^{ - 1}} $.
 		 Moreover, the universal and  implied constants only depend on $T $ and  $ {\max _{0 \leqslant n \leqslant T}}\left| {f\left( n \right)} \right| $\footnote{Here we regard $ p $ and $ q $ as fixed.}. 
 	\end{itemize}
 	\end{itemize}
 \end{theorem}
\begin{remark}
As previously emphasized, although the results here are stated specifically for the case of toral translation, they can be readily formulated in a conjugate manner \cite{DSSY17,DY18,DM23,RB24} to extend to general nonlinear dynamical systems (for example, those that can be smoothly conjugated to toral translations through KAM theory, etc.). Consequently, Theorem \ref{QMC} has broad applicability, particularly in celestial mechanics, since the motion of celestial bodies is often quasi-periodic. For the sake of brevity, we prefer not to complicate the issue further.
\end{remark}

 To facilitate use and avoid verbosity, Theorem \ref{QMC} primarily synthesizes some (not all) important existing results.  It is crucial to emphasize that, according to the results above, the \textit{theoretical} exponential convergence rate of this weighted quasi-Monte Carlo method increases as $ \min\{p, q\} $ grows. However, we cannot be certain that this is always the case in \textit{practical} applications, as we must take into account the inherent computational errors of the computer itself\footnote{On many occasions, data may overflow, leading to reduced precision or even rendering calculations impossible.}.

  The terms ``quasi-periodic'' and ``almost periodic'' are used here because the corresponding rotations are nonresonant; for example, they satisfy certain finite-dimensional or infinite-dimensional Diophantine conditions as in Definitions \ref{Finite-dimensional Diophantine nonresonance} and \ref{Infinite-dimensional Diophantine nonresonance}. To simplify the exposition, we prefer not to state these concepts in Theorem \ref{QMC}. These rotations are universal in the sense of measure, meaning that the set of them admits full Lebesgue measure or full probability measure, as previously mentioned. This weighted quasi-Monte Carlo method is \textit{not} only effective for quasi-periodic, almost periodic, and periodic cases, but also yields similar results for general decaying waves \cite[Theorems 1.1 and 1.2]{TL25b}. As a corollary, when the observational data can be decomposed into a linear combination of quasi-periodic, almost periodic, and periodic components, as well as decaying waves, this method also exhibits excellent acceleration effects. Here, we also mention that the ``super-analyticity'' imposed for the quasi-periodic case in (II)--(6) of Theorem \ref{QMC} is not as strong as it might appear at first glance (of course, the corresponding ``super-analyticity'' in the finite-dimensional case is indeed strong). This is because, unlike the finite-dimensional case, the infinite-dimensional case requires consideration of the cardinality estimates brought by the $\mathbb{Z}_*^{\infty}$ lattice, which inevitably leads to stronger decay of individual Fourier coefficients.
 
  \begin{figure}[h] 	\centering 	\includegraphics[width=350pt]{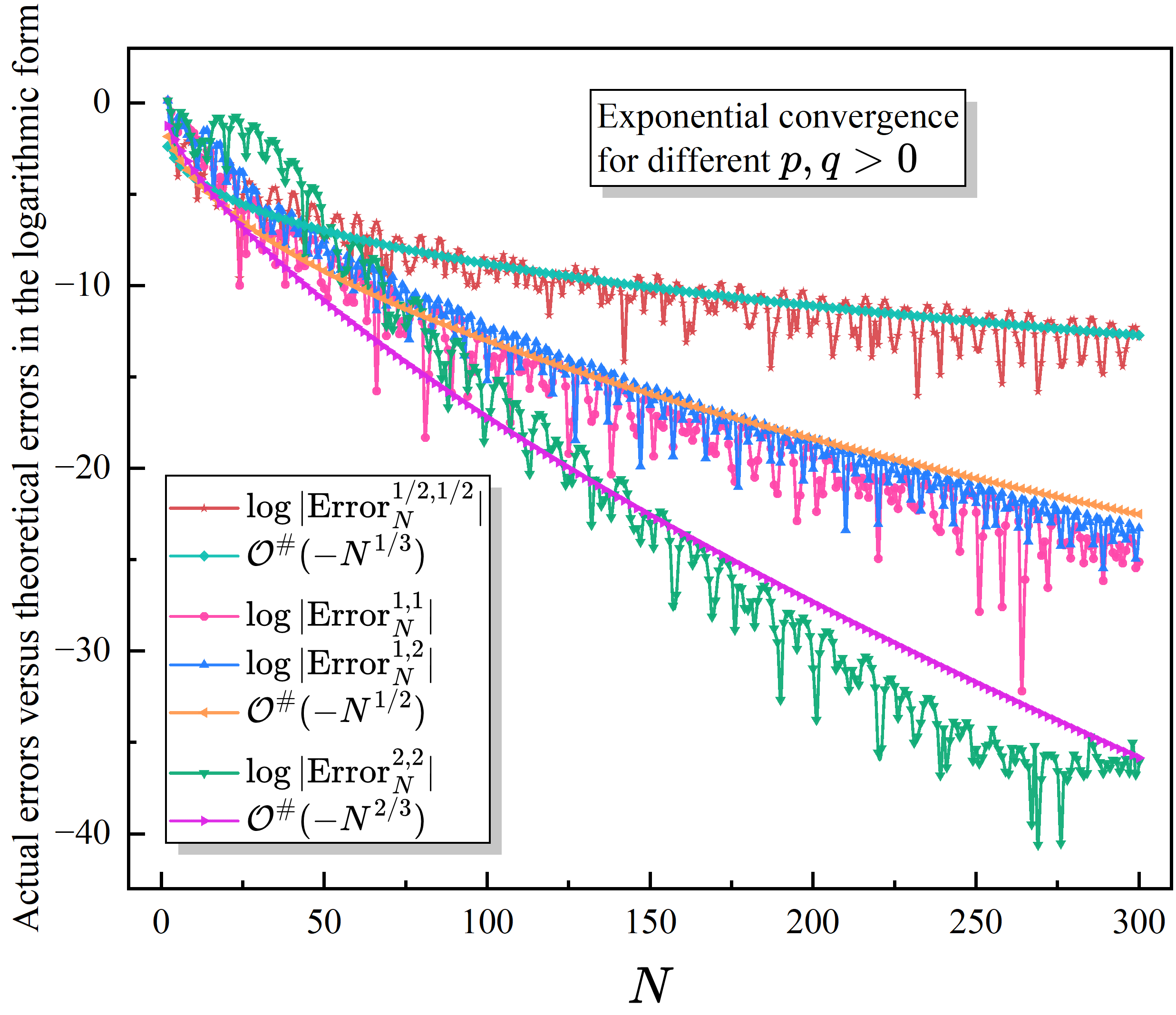} 	\caption {Exponential convergence of the  weighted quasi-Monte Carlo method  in different settings. For the oscillating curves, red corresponds to $ p = q = 1/2 $, pink corresponds to $ p = q = 1 $, blue corresponds to $ p = 1$ and $ q = 2 $, and green corresponds to $ p = q = 2 $, respectively.} 	\label{FIGWBA5}
 	
 \end{figure}
 
 To facilitate the reader's understanding of Theorem \ref{QMC}, we present the numerical simulations in Figure \ref{FIGWBA5}. For convenience, we denote the actual error in the $\ell^{\infty}$-norm as ${\rm{Error}}^{p,q}_N$.
 We consider a quasi-periodic orbit generated by $f(x) = \sin(2\pi x) + \cos(10\pi x)$, with an initial value $ \theta=0 $ and a rotation number  $\rho = 1/2\pi$, which is universal (i.e., not specifically constructed to achieve the ``best'' convergence). The choice of $f$ as a finite-order trigonometric polynomial is motivated by practical considerations: in numerical simulations, we cannot truly achieve infinity, and extremely high orders would lead to inaccurate computations.
 According to Theorem \ref{QMC}, specifically the part (I)--(4), the theoretical error is $\mathcal{O}\left(\exp(-cN^{\zeta_2})\right)$, where $ {\zeta _2} = {\left( {1 + 1/\min \left\{ {p,q} \right\}} \right)^{ - 1}} $ and $c>0$ is a universal constant that may depend on $p$ and $q$. To visualize the results more intuitively, we take the \textit{logarithm} of both the actual error and the theoretical error (and thus the theoretical error becomes $\mathcal{O}\left(-N^{\zeta_2}\right)$, where $ {\zeta _2} = {\left( {1 + 1/\min \left\{ {p,q} \right\}} \right)^{ - 1}} $). This transformation results in curves that are more curved than $\mathcal{O}^{\#}(-N)$. We present four cases: $p = q = 1/2$, $p = q = 1$, $p = 1$ and $q = 2$, and $p = q = 2$. It can be  observed from Figure \ref{FIGWBA5} that as $\min\{p, q\}$ increases, the exponential convergence rate of the actual error indeed accelerates as we proved, closely matching the changes in the order of the theoretical error (although we do not estimate the control constant $ c>0 $ as in \cite{TL25b} for general decaying waves, which would have provided a more precise result). Notably, for the cases $p = q = 1$ and $p = 1$, $q = 2$, the actual errors are \textit{almost identical}, aligning well with our theoretical conclusions, as $\min\{1, 1\}=\min\{1, 2\}=1$ at this point.

 In addition to the results in the general sense mentioned above, \cite[Theorems 2.1, 2.2, 3.1 and 3.2]{TL24a} and \cite[Theorems 3.1--3.4]{TL24b} also investigate results based on Liouville rotations. Specifically, for any Liouville type rotation, as long as the observable is sufficiently smooth, the weighted Birkhoff average can converge at an exponential rate. This is quite different from the phenomena observed in classical dynamical systems: Firstly, many dynamical systems essentially impose restrictions on the irrationality of rotations, such as the Brjuno condition, and \textit{do not} allow arbitrary Liouville rotations. Secondly, even for the specific case of irrational translations, without weighting, such results are absolutely impossible according to the co-homological equation argument. However, when the regularity of the observable is insufficient, such weak nonresonance will reduce the convergence rate in the weighted case, which is entirely reasonable from the perspective of dynamical systems.

 Finally, we would like to announce that this weighted quasi-Monte Carlo method can be further optimized by exploring other weighting functions. For the sake of brevity, we present only one specific case below.
 
 \setcounter{footnote}{0}
 \renewcommand{\thefootnote}{\fnsymbol{footnote}}
 \begin{theorem}\label{TH3}
 	There exist weighting functions $ w^*(x) \in C_0^\infty([0,1]) $ such that the following holds. If the observable $ f $ is sufficiently smooth on $ \mathbb{T}^d $ with $ d \in \mathbb{N}^+ \cup \{  + \infty \}  $, then for almost all rotations $ \rho \in \mathbb{R}^d $, the uniform convergence of the corresponding\footnote{The normalization constant $ A^{p,q}_N $ also changes accordingly.} weighted Birkhoff average \eqref{WBA} exhibits a nearly linear-exponential rate. More precisely, the error is $ \mathcal{O}\left( {\exp ( { -c {N/\log N}} )} \right) $ for some universal constant $ c>0 $. An explicit example of $ w^*(x) $ is given by
\begin{equation}\notag
	{w}^*\left( x \right):= \left\{ \begin{array}{ll}
		\mathscr{A} \exp \left( { - \exp ( {{x^{ - 1}}} ) - \exp ( {{{( {1 - x} )}^{ - 1}}} )} \right),&x \in \left( {0,1} \right),  \hfill \\
		0,&x \notin \left(0,1\right), \hfill \\
	\end{array}  \right.
\end{equation}
where the normalization constant is $ \mathscr{A} := {\left( {\int_0^1 { \exp \left( { - \exp ( {{s^{ - 1}}} ) - \exp ( {{{\left( {1 - s} \right)}^{ - 1}}} )} \right)ds} } \right)^{ - 1}}.$
 \end{theorem}
\begin{remark}\label{RE13}
Here, the logarithmic part in $ \mathcal{O}\left( {\exp ( { -c {N/\log N}} )} \right) $ could be replaced by any approximation function that is monotonically increasing and tends to $ +\infty $ (as $ x \to +\infty$), by selecting other weighting functions. However, we conjecture that the linear-exponential rate $ \mathcal{O}\left( {\exp ( { -c {N}} )} \right) $ cannot be achieved.
\end{remark}
 Taking into account (I)--(3) in Theorem \ref{QMC}, the conclusion of Theorem \ref{TH3} is, in fact, \textit{not} counterintuitive. The proof of this theorem essentially relies on the novel techniques established in \cite[Lemma 4.1]{TL25b}, and the overall proof is similar to the quantitative results presented in \cite[Theorem 6.1]{TL25b}. However, it should be emphasized that no other existing techniques can effectively prove Theorem \ref{TH3} beyond this approach. We have decided to present the complete proof in a separate paper and will simultaneously  investigate the lower bounds of the convergence rates for the corresponding weighted Birkhoff averages.

\subsubsection{Physically related problems  revisited}\label{SECphy}
 Recall that Theorem \ref{QMC}  provides the theoretical convergence rates of the weighted quasi-Monte Carlo method utilized in \eqref{WBA}. To better illustrate its  practical effectiveness in deterministic dynamical systems, we present in this section a selection of specific problems related to physics. These problems, which have recently been studied in detail through illustration and numerical simulation, encompass a wide array of highly significant research subjects, including the (circular/planar) restricted three-body problem, the Arnold circle map, the standard map, the Lorenz map, the H\'enon map, the Duffing oscillator, and many more. For the sake of simplicity, we will assume that the time scale is $ N $ in the following discussions.

\begin{itemize}
\item[(I)] Computing the integral of a periodic smooth mapping \cite{DY18}: In the periodic case, one can artificially construct quasi-periodicity. This coincides with the concept of the weighted quasi-Monte Carlo method. For periodic mappings with a finite number of variables, one can normalize them to the standard torus via scaling. Then, by selecting quasi-periodic translations with Diophantine rotations and performing the weighted quasi-Monte Carlo method, one can efficiently compute the integral of periodic mappings, if the mappings are smooth (e.g., $ C^\ell $, $ C^\infty $ or even analytic). It is also worth noting that the irrationality of the rotations can be  specified \textit{arbitrarily}. Therefore, it is preferable to construct rotations with \textit{critical} Diophantine exponents. These special rotations, form a set which admits zero Lebesgue measure but full Hausdorff dimension, \textit{maximize} computational efficiency.
 
\item [(II)] Computing rotations \cite{DSSY17,DY18,SY18,SM20,MS21,BC23,DM23,BM24,CCGd24,BdJC25,BHS25,MS25,SM25a,SM25b}: This is likely the \textit{most} extensively utilized application of this weighted quasi-Monte Carlo method currently. The calculation of rotations---which are always well-defined and independent of the initial value---is often \textit{not} an isolated issue. Instead, it may also serve as an \textit{essential} intermediate component in other computational processes. High-precision rotations calculated through exponential acceleration enable one to often rely on \textit{experience} to determine their arithmetic properties, thereby judging the corresponding dynamic properties. In the unweighted case, the convergence rate for computing rotations using standard Birkhoff averages is typically at most $ \mathcal{O}\left(N^{-1}\right) $, which stands in stark contrast to the exponential convergence achieved with weighting. The extensive use of this weighted quasi-Monte Carlo method in numerous studies strongly attests to its effectiveness.

\item [(III)] Computing Lyapunov exponents \cite{DSSY17,SM20,DM23,MS25,SM25a}:  In dynamical systems, positive Lyapunov exponents, which measure average exponential stretching along orbits, typically characterize chaos. In contrast, for quasi-periodic/almost periodic systems, the Lyapunov exponents are zero. If all Lyapunov exponents are strictly less than zero, the system is certainly not chaotic but rather asymptotically stable.  However, the standard orthogonalization-based method for computing Lyapunov exponents converges very slowly---typically at most 
$ \mathcal{O}\left(N^{-1}\log N\right) $---and may not converge at all. This weighted quasi-Monte Carlo method, while not guaranteed, can \textit{sometimes} significantly accelerate the computation of Lyapunov exponents. However, there is currently no rigorous proof for this phenomenon, possibly due to singularities from the logarithmic term. This is an important issue we will need to address in the future.

\item [(IV)] Distinguishing between regular and chaotic orbits \cite{SM20,MS21,DM23,BHS25,BBC25,GZE25,MS25,SM25b,ZXH25}: This is one of the most fundamental questions in the field of dynamical systems\footnote{Prof. James D. Meiss's research group has made many fundamental contributions to (III) and (IV).}. Although the Lyapunov exponent can be used to determine whether a system is chaotic to some extent, it is \textit{not} a necessity. There are many other techniques that can achieve this goal, and most of them do not require the computation of specific Lyapunov exponents. This weighted quasi-Monte Carlo method can \textit{sometimes} effectively distinguish between (weak/strong) chaotic orbits, island chains, and quasi-periodic/almost periodic orbits. For chaotic cases, in general, its convergence rate is no faster than that of the unweighted average of a random signal, with an error of $ \mathcal{O}\left(N^{-1/2}\right) $. In contrast, for other regular cases, the convergence rate is typically exponential,  provided the system is analytic, as detailed in Theorem \ref{QMC}.

\item [(V)] Computing Fourier coefficients \cite{DSSY17,DY18,SY18,BM24,Boh24,RB24,RKB25}: In addition to its extremely broad practical applications, it also involves an old problem of determining the regularity of the  conjugacy: if the Fourier coefficients decay exponentially, then one \textit{intuitively} concludes that the conjugacy is analytic. In the subsequent   section, we establish   for the first time a rigorous result (Theorem \ref{CLT-FOU}) that demonstrates the following: if the underlying system is analytic and the rotation is Diophantine, then this weighted quasi-Monte Carlo method can indeed accelerate the computation of Fourier coefficients of the embedding exponentially.

\end{itemize}

In addition to the aforementioned applications, this weighted quasi-Monte Carlo method can also be employed to numerically validate the error of the ``label functions'' \cite{RB25}. It has also been extended to other fields, including the study of numerical methods and other consistent approximations of the Euler system in gas dynamics, via a new concept of (S)-convergence \cite{Fei21}.  We would also like to mention \cite{RB24}, which employs reduced rank extrapolation to explore more general weighting methods (not necessarily in the form of \eqref{WBA}, but depending on the specific systems), and achieves arbitrary exponential convergence similar to (I)--(3) in Theorem \ref{QMC}. In a study related to that of \cite{RB24} but developed independently, \cite{Sal24} introduced graph filters from the perspective of graph signal processing theory and explored the fast convergence in the ergodic theorem for Markov chains. Recently, \cite{BBC25} presented a comprehensive integration of this weighted quasi-Monte Carlo method with five data-driven algorithms from dynamical systems, accompanied by extensive examples.  We are convinced that this weighted quasi-Monte Carlo method holds potential applications far beyond what has been discussed in this paper. We look forward to seeing more works that delve into a broader range of problems related to dynamical systems to further demonstrate its effectiveness.

\subsubsection{Exponentially computing Fourier coefficients using weighted Birkhoff averages}\label{SECFOU}
To facilitate the reader's understanding, we retain the primary notations from \cite{BM24} (pertaining to the one-dimensional torus) and further explore the more general cases on higher-dimensional and infinite-dimensional tori.

Assume that $ d, D \in  \mathbb{N}^+ \cup \{  + \infty \} $. For a diffeomorphism $ F:\Omega  \to \Omega  $ with $ \Omega  \subset {\mathbb{R}^D} $, suppose that it admits an \textit{analytic}\footnote{Recall from Definitions \ref{Infinite-dimensional analyticity} and \ref{Infinite-dimensional Diophantine nonresonance} the concepts of analyticity in the finite-dimensional and infinite-dimensional settings, respectively.} quasi-periodic/almost periodic invariant curve $ \Gamma $. An important issue is to compute the Fourier coefficients of a lift/parameterization for $ \Gamma $. Under these grounds, there exists an \textit{analytic} mapping $ K:{\mathbb{T}^d} \to {\mathbb{R}^D} $ such that $ {\rm{image}}{K}=\Gamma $, with the following conjugation:
\begin{equation}\label{CONJU}
	F\left( {K\left( \theta  \right)} \right) = K\left( {\theta  + \rho } \right), \quad  \text{where $ \rho \in \mathbb{T}^d $ is a nonresonant rotation.}
\end{equation}
Here, we denote $ K(\theta_0):=p_0 $ for some initial angular value $ \theta_0 $.  Typically, $ \rho $ can also be computed with exponential convergence using the weighted quasi-Monte Carlo method, as explained in Section \ref{SECphy}. Figure \ref{FIG3} illustrates the geometric intuition of the conjugation \eqref{CONJU} for the case $ d = 1 $, as well as the algebraic commutative diagram for the general case $ d \in \mathbb{N}^+ $.

\begin{figure}[htbp]
	\centering 
	\begin{minipage}[c]{0.49\textwidth} 
		\centering 
		\begin{tikzpicture}[scale=1.4] 
		\draw[->,  thick, -{Latex[scale=1.1]}] (-1cm,0cm) -- (5cm,0cm) node[right]{\(x\)}; 
		\draw[->,  thick, -{Latex[scale=1.1]}] (0cm,-1cm) -- (0cm,4.5cm) node[above]{\(y\)}; 
		
		\fill (0cm,0cm) circle (1.5pt) node[below left]{\(O\)};

		\draw[fill=red!15, semithick, decoration={markings,
			mark=between positions 0.1 and 0.9 step 0.1 with {\arrow{latex}},
		},postaction={decorate}
		] plot[smooth, tension=0.5] coordinates {
			(3cm+0.3cm,2.4cm+0.3cm) 
			(2.5cm+0.3cm,3cm+0.3cm) 
			(1.7cm+0.3cm,2.7cm+0.3cm) 
			(1cm+0.3cm,3cm+0.3cm) 
			(0.5cm+0.3cm,2.5cm+0.3cm) 
			(0.3cm+0.3cm,1.5cm+0.3cm) 
			(1cm+0.3cm,1.3cm+0.3cm) 
			(1.5cm+0.3cm,0.5cm+0.3cm) 
			(2.5cm+0.3cm,0.6cm+0.3cm) 
			(3cm+0.3cm,0.8cm+0.3cm) 
			(3.8cm+0.3cm,1.5cm+0.3cm) 
			(3.5cm+0.3cm,2.5cm+0.3cm)
			(3.3cm+0.3cm,2.6cm+0.3cm)
		} -- cycle;
		
		
		\fill[purple!90] (3cm+0.3cm,2.4cm+0.3cm) circle (1.2pt) node[above right, yshift=0.2cm, xshift=-0.2cm, black] {\(K(\theta_0):=p_0\)};
		
		\fill[mydeepgreen]  (1cm+0.3cm,3cm+0.3cm) circle (1.2pt) node[above,black] {\(K(\theta_0+\rho):=p_1\)};
		
		\draw[dashed, ->, blue, semithick] (2cm+0.3cm,1.8cm+0.3cm) -- (4.0cm+0.3cm,1.8cm+0.3cm) node[below] {\(\theta\)};
		
		\draw[dashed, ->, blue, semithick] (2cm+0.3cm,1.8cm+0.3cm) -- (2cm+0.3cm,3.8cm+0.3cm);
		
		\draw[dashed, <-, blue, semithick] (0cm+0.3cm,1.8cm+0.3cm) -- (2.0cm+0.3cm,1.8cm+0.3cm);
		
		\draw[dashed, ->, blue, semithick] (2cm+0.3cm,1.8cm+0.3cm) -- (2cm+0.3cm,-0.2cm+0.3cm);
		
		\draw[-latex,purple!90,semithick] (2cm+0.3cm,1.8cm+0.3cm) -- (3cm+0.3cm,2.4cm+0.3cm);
		
		\draw[-latex,mydeepgreen,semithick] (2cm+0.3cm,1.8cm+0.3cm) -- (1cm+0.3cm,3cm+0.3cm);
		
		\tikzset{
			arc arrow/.style = {
				decoration = {
					markings,
					mark = at position 0.68 with {\arrow{latex}}
				},
				postaction = decorate
			}
		}
		
		\draw[purple!90, semithick, arc arrow] (2.9cm+0.3cm,1.8cm+0.3cm) arc (0:40:0.69cm)
		node[midway, below right, yshift=0.25cm] {\scalebox{0.8}{\(\theta_0\)}};

		\tikzset{arc arrow/.style = {
				decoration = {markings,
					mark = at position 0.55 with {\arrow{latex}}
				},
				postaction = decorate
		}}
		
		\draw[mydeepgreen, semithick, arc arrow] (1.8cm+0.3cm,1.8cm+0.3cm) ++(0:1cm) arc (10:115:1cm)
		node[midway, below left, xshift=0.25cm] {\scalebox{0.8}{\(\theta_0+\rho\)}};
		

		\fill (2cm+0.3cm,1.8cm+0.3cm) circle (1.5pt) node[below left] {$O'$};
		\end{tikzpicture}
	\end{minipage}
	\hfill 
	\begin{minipage}[c]{0.49\textwidth} 
		\centering 
		\begin{tikzcd}[row sep=3.9 cm, column sep=3.9 cm, arrows={thick}]
			{\theta}  \arrow[r, " \mathcal{T}\; {\text{(with rotation $ \rho $)} }"{font=\large}] \arrow[d, "K"'{font=\large}] & \theta+\rho \arrow[d, "K"{font=\large}] \\
			K(\theta) \arrow[r, "F"'{font=\large}] & K(\theta+\rho)
		\end{tikzcd}
	\end{minipage}
	\caption{The geometric intuition of the conjugation \eqref{CONJU} (with $ d=1 $) and the commutative diagram}
	\label{FIG3}
\end{figure}

 Due to the analyticity of $ K $, we can represent it as a Fourier expansion (which is also a vector of $ D $ components):
\[K\left( x \right) = {\left( {{K^1}\left( x \right), \ldots ,{K^D}\left( x \right)} \right)^{\top}} = \sum\limits_{k \in {\mathbb{Z}^d}} {\widehat K\left( k \right)\exp \left( {2\pi i\left\langle {k,x} \right\rangle } \right)}  \quad \text{or} \quad  \sum\limits_{k \in \mathbb{Z}_ * ^\infty } {\widehat K\left( k \right)\exp \left( {2\pi i\left\langle {k,x} \right\rangle } \right)} \]
with 
\[\widehat K\left( s \right) = \int_{{\mathbb{T}^d}} {K\left( x \right)\exp \left( { - 2\pi i\left\langle {s,x} \right\rangle } \right)dx} .\]
Here, for the infinite-dimensional case with an observable $ u $ possessing at least analyticity, 
\[\int_{{\mathbb{T}^\infty }} {u\left( x \right)dx} : = \mathop {\lim }\limits_{\mathcal{L} \to  + \infty } \int_{{\mathbb{T}^\mathcal{L}}} {u\left( x \right)d{x_0} \ldots d{x_\mathcal{L}}} \]
is well-defined and equals $ \widehat u(0) $, see \cite{Bou05,MP21,CGP24}, for instance.
 Set $ {\theta _n}: = {\theta _0} + n\rho  $ for $ n \in \mathbb{N} $. Then, by induction, we obtain from the conjugation \eqref{CONJU} that 
\begin{equation}\label{PNDY}
	{p_n}: = {F^n}\left( {{p_0}} \right) = {F^n}\left( {K\left( {{\theta _0}} \right)} \right) = {F^{n - 1}}\left( {K\left( {{\theta _0} + \rho } \right)} \right) =  \cdots  = K\left( {{\theta _0} + n\rho } \right) = K\left( {{\theta _n}} \right).
\end{equation}
As an application of the weighted quasi-Monte Carlo method, we claim that for fixed $ s \in \mathbb{Z}^d $ (or $ \mathbb{Z}_*^\infty $),
\begin{equation}\label{ill}
	\widehat K\left( s \right) \approx \exp \left( { - 2\pi i\left\langle {s,{\theta _0}} \right\rangle } \right) \cdot \frac{1}{{A_N^{p,q}}}\sum\limits_{n = 1}^{N - 1} {{w_{p,q}}\left( {n/N} \right){p_n}\exp \left( { - 2\pi in\left\langle {s,\rho } \right\rangle } \right)},  \quad  N \gg 1,
\end{equation}
with \textit{exponential} convergence. However, there is currently no rigorous proof regarding its convergence rate, and the formal derivation of the above formula is not sufficiently rigorous. To refine the theory in this aspect, we establish Theorem \ref{CLT-FOU} as follows, providing a nearly  thorough result for \eqref{ill} in the \textit{analytic} setting.

\begin{theorem}\label{CLT-FOU}
	For almost all rotations $ \rho \in \mathbb{T}^d $ with $  d \in  \mathbb{N}^+ \cup \{  + \infty \} $, we have that for fixed $ s \in \mathbb{Z}^d $,
\begin{equation}\label{FOUSL}
	{\left\| {\widehat K\left( s \right) - \exp \left( { - 2\pi i\left\langle {s,{\theta _0}} \right\rangle } \right) \cdot \frac{1}{{A_N^{p,q}}}\sum\limits_{n = 1}^{N - 1} {{w_{p,q}}\left( {n/N} \right){p_n}\exp \left( { - 2\pi in\left\langle {s,\rho } \right\rangle } \right)} } \right\|_{{\ell ^\infty }}} = o\left( 1 \right)
\end{equation}
 holds uniformly in $ \theta_0 \in \mathbb{T}^d $ as $ N \to +\infty $,	with $ p_n $ defined in \eqref{PNDY}. Moreover, we have: 
\begin{itemize}
\item[(I)] The finite-dimensional case: the error in \eqref{FOUSL} is $ \mathcal{O}\left( {\exp ( - {N^{{\zeta _6}}})} \right) $ for any 
\[0 < {\zeta _6} < {\left( {d + 1 + 1/\min \left\{ {p,q} \right\}} \right)^{ - 1}}.\]
 In particular, it holds uniformly in $ \|s\|_{\ell^1}=o\left({N^{{\zeta _6}}}\right)  $.

\item[(II)] The infinite-dimensional case: the error in \eqref{FOUSL} is $ \mathcal{O}\left( {\exp ( - {{(\log N)}^{{\zeta _7}}})} \right) $ for any  $ 2 \leqslant {\zeta _7} < 1 + \eta $. In particular, it  holds uniformly in $ |s|_\eta=o\left({( {\log N} )^{\zeta_7} }\right)  $.
\end{itemize}
\end{theorem}
\begin{remark}\label{RE14}
We emphasize that Theorem \ref{CLT-FOU} is by no means a direct corollary of Theorem \ref{QMC}. On the contrary, portions of the proof of Theorem \ref{QMC} are embedded within the more complicated argument required for Theorem \ref{CLT-FOU}. A detailed discussion will be given in Section \ref{SECPT12}.
\end{remark}
The conclusion and proof of Theorem \ref{CLT-FOU} are distinct from existing results. Moreover, it determines the \textit{effective order} of $ s $ (with respect to the time scale $ N $) in computing Fourier coefficients, under the premise of ensuring exponential convergence. It is worth emphasizing that this is highly crucial in practical numerical simulations \cite{DSSY17,SM20,DM23,MS25,SM25a}, as it indicates how many orders of Fourier coefficients can be rapidly and accurately calculated within the specified time scale. Certainly, the convergence rate can still be estimated for $ s $ beyond the effective order, but perhaps it will be relatively slower (e.g., of the polynomial type), which is consistent with both intuition and numerical simulations. One can also utilize Theorem \ref{QMC} to discuss more general cases, which we do not pursue in this paper.

\subsubsection{Which factors may affect the ``practical'' convergence rate?}\label{SEC114}

In contrast to the dynamical-systems perspective, researchers  focusing on numerical simulations are typically more interested in the ``practical'' convergence rate for finite yet large $N$. This motivates the following question: \textit{Which  factors may affect the ``practical'' convergence rate of the weighted quasi-Monte Carlo method investigated in this paper?} While \cite[Sections 4 and 5]{DSSY17} offer  relevant insights, we provide in this section a deeper and more theoretical discussion.

Let us consider the weighted Birkhoff average in \eqref{WBA} utilizing the  weighted quasi-Monte Carlo method with $ p,q>0 $. As in Theorem \ref{QMC}, we continue to examine its uniform convergence rate to $\int_{\mathbb{T}^d} f   dx$ in the $\ell^\infty$-norm.  In order to clarify the presentation, we focus exclusively on the case $d\in \mathbb{N}^+$, and our convergence analysis is restricted to the polynomial type of finite order\footnote{Given that this paper aims to make this weighted quasi-Monte Carlo method more user-friendly, a more precise quantitative analysis of exponential convergence (from a more dynamical-systems perspective) will be provided as a byproduct in a separate paper focusing on large deviations.}. According to \cite{DSSY17,DY18,DM23,TL24a}, the convergence rate is   $\mathcal{O}\left(N^{-m}\right) $ when the regularity of the observable $f$ and the nonresonance of the rotation $\rho$  satisfies
\[\sum\limits_{0 \ne k \in {\mathbb{Z}^d}} {{{ \| {\widehat f(k)}  \|}_{{\ell ^\infty }}}\left( {\mathop {\sup }\limits_{n \in \mathbb{Z}} \left( {{{\left| {\left\langle {k,\rho } \right\rangle  - n} \right|}^{ - m}}} \right)} \right)}  <  + \infty  ,\quad  2 \leqslant m \in \mathbb{N}^+. \]
Observe that if $\rho$ is nonresonant, we can assume without loss of generality that it satisfies the following nonresonance condition 
\[\left| {\left\langle {k,\rho } \right\rangle  - n} \right| > \frac{\gamma }{{\Delta \left( {{{\left\| k \right\|}_{{\ell ^1}}}} \right)}},\quad \gamma  > 0,\quad \forall 0 \ne k \in {\mathbb{Z}^d},\quad \forall n \in \mathbb{Z},\]
where $\Delta$ is an approximating function (see Remark \ref{RE13}) with $\Delta(1) > 0$. In this case, using a slightly modified argument in \cite{DSSY17,DY18} (for instance, \cite[Section 4.1]{TL24a}), the convergence rate can be estimated as $ \big( {\prod\nolimits_{j = 1}^4 {{{\mathcal  C}_j}} }  \big)N^{-m}$ ($2 \leqslant m \in \mathbb{N}^+$), where the coefficients are explicitly given by
\[{\mathcal{C}_1} = \mathop {\sup }\limits_{N \in {\mathbb{N}^ + }} \frac{N}{{A_N^{p,q}}},\quad {\mathcal{C}_2} = {\left\| {{D^m}{w_{p,q}}} \right\|_{{L^1}\left( {0,1} \right)}},\quad {\mathcal{C}_3} = \frac{1}{{{{\left( {2\pi \gamma } \right)}^m}}},\quad {\mathcal{C}_4} = \sum\limits_{0 \ne k \in {\mathbb{Z}^d}} {{{ \| {\widehat f(k)}  \|}_{{\ell ^\infty }}}{\Delta ^m}\left( {{{\left\| k \right\|}_{{\ell ^1}}}} \right)} .\]
At this point, we can identify  which factors may affect the convergence rate. Please keep in mind that this serves only as an \textit{intuitive} understanding and a \textit{possible} explanation for numerical simulations, rather than a rigorous mathematical proof---because the convergence rate discussed here is not necessarily the actual rate, but merely an upper bound. To study this problem rigorously, one must simultaneously consider the lower bound of the convergence rate, which requires entirely different techniques and is considerably more challenging. This will be detailed in \cite{TL25c}. However, from the perspective of numerical simulations, the following interpretation from a dynamical-systems perspective may be helpful. Indeed, the following analysis has been empirically validated as correct in practice.
\begin{itemize}
	\item[(I)] In numerical simulations, the convergence rate is sometimes approximated by finite-order polynomials such as $N^{-m}$ (even though it might actually be exponential) for finite yet large $N $. Hence, the specific value of $m$ is a potential factor.
	
	\item[(II)] The coefficients $ {\mathcal{C}_1} $ and $ {\mathcal{C}_2} $ reveal  that the weighting function $w_{p,q}(x)$ itself is a potential   factor. Straightforward verification shows  $\mathop {\lim }\nolimits_{N \to  + \infty } N/A_N^{p,q} = 1$ (see Section \ref{SEC14}), but for large $p$ and $q$, finite $N$ may cause $ {\mathcal{C}_1} $ to deviate significantly from $1$. As for $ {\mathcal{C}_2} $, \cite[Lemma 4.1]{TL25b} provides its asymptotic estimate with respect to $m$, allowing for a clearer assessment of its potential impact.
	
	\item[(III)]  The coefficients $ {\mathcal{C}_3} $ and $ {\mathcal{C}_4} $ reveal  that the rotation  $\rho$ is also a potential factor. This can be explained from two aspects. First, for a given nonresonance condition---that is, with $\Delta$ fixed---a smaller $\gamma$ leads to a larger $ {\mathcal{C}_3} $. Second, for a fixed $\gamma$, the type of nonresonance---such as Diophantine, Brjuno, or Liouville---also affects $ {\mathcal{C}_4} $. In particular, if we further consider a given observable $f$, then from $ {\mathcal{C}_4} $ we can intuitively  see that weaker nonresonance tends to slow down the ``practical'' convergence rate. For example, when $\rho$ is of Diophantine type, i.e., $\Delta (x)= x^\tau$ with some $\tau \geqslant d$, then $C_4$ becomes $ \sum\nolimits_{0 \ne k \in {\mathbb{Z}^d}} {{{ \| {\widehat f(k)}  \|}_{{\ell ^\infty }}}\left\| k \right\|_{{\ell ^1}}^\tau } $. Clearly, sufficiently fast polynomial decay of the Fourier coefficients of $f$ ensures that $C_4 < +\infty$. However, if $\rho$ is of Liouville type, such as when $\Delta (x)= \exp(x)$ (and  there simultaneously exist infinitely many tuples $ \left( {k,n} \right) \in {\mathbb{Z}^d} \times \mathbb{Z} $ such that $ \left| {\left\langle {k,\rho } \right\rangle  - n} \right| \sim \gamma \exp \left( { - {{\left\| k \right\|}_{{\ell ^1}}}} \right) $), then to guarantee $C_4 =\sum\nolimits_{0 \ne k \in {\mathbb{Z}^d}} {{{ \| {\widehat f(k)}  \|}_{{\ell ^\infty }}}\exp \left( {{{\left\| k \right\|}_{{\ell ^1}}}} \right)} < +\infty$, we have to  require exponential decay of the Fourier coefficients of $f$. Theoretical results concerning Liouville nonresonance have been extensively discussed in \cite{TL24a,TL24b,TL25a}.  \cite{DSSY17} also corroborates these viewpoints through numerical simulations.   As a remark, in practice there appears to be no feasible way to decide whether $\rho$ is Liouville; one usually has to rely on empirical judgement\footnote{Z. Tong first came across this point during a 2024 talk by Prof. Alessandra Celletti related to weighted Birkhoff averages.}.   Consequently, both aspects above are indispensable.  For instance, even if $\rho$ is Diophantine, an extremely small $\gamma$ may still tempt one to place it, intuitively, in the Liouville class---although, from a dynamical-systems perspective, the two situations are entirely different.

	\item[(IV)] The coefficient $ {\mathcal{C}_4} $ reveals  that $f$ is also a potential factor. This can be examined from two perspectives.  
	On the one hand, given a fixed rotation $\rho$, it is clear that higher regularity of $f$ increases the possibility of $  {\mathcal{C}_4}   < +\infty$. However, this may not significantly affect the size of $ {\mathcal{C}_4} $, since regularity primarily affects the decay of Fourier coefficients at high frequencies, while the value of $ {\mathcal{C}_4} $ is partly determined by low-frequency Fourier coefficients.    On the other hand, the only rigorous statement that can be made here concerns how the   ``scale'' of $f$   affects the ``practical'' convergence rate---for example, replacing $f$ with $2f$. The homogeneity of the weighted quasi-Monte Carlo method with respect to $f$ provides a direct explanation in such cases.
\end{itemize}

In conclusion, the ``practical'' convergence rate may be affected by the specific value of $m$, as well as all dynamical system-related factors involved in this weighted quasi-Monte Carlo method---namely, the weighting function, the rotation, and the observable. It can be seen that the most critical factors are the \textit{nonresonance} of the rotation and the \textit{regularity} of the observable, which precisely correspond to the central perspectives of interest in dynamical systems.

\subsection{The probabilistic perspective}\label{SECPRO}
Given that this weighted quasi-Monte Carlo method can have an unexpected acceleration effect in some deterministic cases,  a natural question arises as to whether a similar phenomenon occurs in the random setting via the weighting function $ 	{w}_{p,q}\left( x \right) $ defined in  \eqref{CLTFUN}. Fortunately, there has been some relevant work in history addressing this issue.

\subsubsection{Weighted laws of large numbers and weighted strong laws of large numbers}\label{SECWL}

The weighted laws of large numbers and the weighted strong laws of large numbers have long been central topics in probability theory. A substantial body of historical and recent advances exists on these topics; see, for instance, \cite{JOP65,Pru66,Roh71,Thr87,BC00,Jaj03,CG07,JL08,LM10,CS16,FMZ17,NMSA20,Ad24} and the references therein.

Firstly, concerning the weighted versions of the laws of large numbers, such as the Bernoulli type, Khinchin type, Chebyshev type, Markov type, Poisson type, and Borel type (via the weighting function $ w_{p,q}(x) $, and similarly for the following), it is quite straightforward to demonstrate that the results in these weighted cases do not significantly differ from those in the unweighted cases. Moreover, the rate of convergence can be readily estimated, which is historically well-known. Therefore, to avoid verbosity, we prefer not to list all the conclusions explicitly here. We only show the Khinchin type and Chebyshev type results without their proof.

\begin{theorem}\label{SJT1}
	\begin{itemize}
	\item[(I)] (Khinchin) Let  $ X_1,  \ldots  , X_N $ be independent identically distributed random variables. Assume further that $ {\bf{E}}X_j $ are uniformly bounded for all $ 1\leqslant j \leqslant N $. Then for any $ \varepsilon>0 $,
	\[\mathop {\lim }\limits_{N \to  + \infty } {\bf P}\left( {\left| {\frac{1}{{{{A_N^{p,q}}}}}\sum\limits_{n = 1}^N {{w_{p,q}}\left( {n/N} \right){X_n}}  - \frac{1}{A_N^{p,q}}\sum\limits_{n = 1}^N {{w_{p,q}}\left( {n/N} \right){\bf E}{X_n}} } \right| < \varepsilon } \right) = 1.\]
	
	\item[(II)]  (Chebyshev) 	Let  $ X_1,  \ldots  , X_N $ be  independent random variables. Assume further that $ {\bf{E}}X_j $ and $ {\bf{E}}X_j^2 $ are uniformly bounded for all $ 1\leqslant j \leqslant N $. Then for any $ \varepsilon>0 $,
	\[\mathop {\lim }\limits_{N \to  + \infty } {\bf P}\left( {\left| {\frac{1}{A_N^{p,q}}\sum\limits_{n = 1}^N {{w_{p,q}}\left( {n/N} \right){X_n}}  - \frac{1}{A_N^{p,q}}\sum\limits_{n = 1}^N {{w_{p,q}}\left( {n/N} \right){\bf E}{X_n}} } \right| < \varepsilon } \right) = 1.\]
\end{itemize}

	\end{theorem}

Secondly, it should be emphasized that the following weighted strong laws of large numbers  (Theorems \ref{LIMT1} and \ref{LIMT2} below) and the weighted central limit theorem (Theorem \ref{CLTT1} in the next section), are primarily built upon  historical results \cite{Thr87,BC00,CG07,Kla09} and have the potential for further refinement---for instance, in conjunction with the previously mentioned works. In particular, these results can be further explored beyond the independent setting, allowing for \textit{dependent} stochastic processes that will cover a broad range of nonlinear phenomena \cite{FMZ17}.

We first present  a  Marcinkiewcz-Zygmund type strong law of large numbers, which utilizes the results in  \cite{Thr87}.
\begin{theorem}\label{LIMT1}
Let  $ X_1,  \ldots  , X_N $ be  independent identically distributed random variables.	Assume further that $ {\bf{E}}X_1=0 $ and ${\bf{E}}{\left| X_1 \right|^\mu} <  + \infty  $	for some $ \mu >0 $. Then for 
	\begin{equation}\notag
	\sigma  = \left\{ \begin{aligned}
		&	1/\mu  + 1/2,&\mu  \geqslant 2, \hfill \\
		&	1/\mu  + 1,&0 < \mu  < 2,  \hfill \\
	\end{aligned}  \right.
\end{equation}
 there holds
	\[{\bf{P}}\left( {\mathop {\lim }\limits_{N \to  + \infty }  {\frac{1}{{{N^{ \sigma }}}}\sum\limits_{n = 1}^N {\sqrt {{w_{p,q}}\left( {n/N} \right)} {X_n}} } =0} \right) = 1.\]
	
\end{theorem}

The following strong law of large numbers with a single logarithm depends on the result in 	\cite{CG07}, which generalizes the findings in \cite{BC00}.

\begin{theorem}\label{LIMT2}
	Let  $ X_1,  \ldots  , X_N $ be  independent identically distributed random variables. Assume further that $ {\bf{E}}X_1=0 $ and $ {\bf{E}}\left( {{{\left| X_1 \right|}^\mu }{{\left( {\log \left| X_1 \right|} \right)}^{ - \frac{\mu }{2}}}} \right) <  + \infty  $ for some $ \mu>2 $. Then  there exists some $ C_1>0 $, such that 
	\[{\bf{P}}\left( {\mathop { {\lim \sup  } }\limits_{N \to +\infty } \left| {\frac{1}{{   \sqrt {N\log N} }}\sum\limits_{n = 1}^N {\sqrt {w_{p,q}\left( {n/N} \right)} {X_n}} } \right| \leqslant C_1{\sqrt {{\bf{E}}{X_1^2}} }} \right) = 1.\]
\end{theorem}

\subsubsection{Weighted central limit theorem}\label{SECWCLT}
 Below,  we present a weighted central limit theorem using the weighting function $ 	{w}_{p,q}\left( x \right) $ in \eqref{CLTFUN}, building upon the foundational work  \cite{Kla09}.

\begin{theorem}\label{CLTT1}
	Assume that $ X_1,  \ldots  , X_N $ are random variables satisfying $ {\bf{E}}X_j^2=1 $ for $ 1\leqslant j \leqslant N $, and the random vector $ X=\left( X_1,  \ldots  , X_N \right) $ is distributed according to a log-concave and unconditional density $ f : \mathbb{R}^N \to \mathbb{R}_{\geqslant 0}$, i.e., $   \log f $ is concave, and for any $ \left( x_1,  \ldots  , x_N \right)\in \mathbb{R}^N$ and a sign vector $ \left( e_1,  \ldots  , e_N \right) \in \left\{ \pm 1  \right\}^N $, 
	\[f\left( {{x_1}, \ldots ,{x_N}} \right) = f\left( {{e _1}{x_1}, \ldots ,{e _N}{x_N}} \right).\]
	Then, it holds with a universal constant $ C_2>0 $ (may depend on $ p$ and $q $) that 
	\[\mathop {\sup }\limits_{\alpha ,\beta  \in \mathbb{R},\;\alpha  < \beta } \left| {{\bf{P}}\left( {\alpha  \leqslant \frac{1}{{\sqrt {{A^{p,q}_N}} }}\sum\limits_{n = 1}^N {\sqrt {w_{p,q}\left( {n/N} \right)} {X_n}}  \leqslant \beta } \right) - \frac{1}{{\sqrt {2\pi } }}\int_\alpha ^\beta  {{\exp\left(- t^2/2\right)}dt} } \right| \leqslant \frac{{C_2 }}{N}.\]
	Moreover, $ {\sup _{0 < p,q \leqslant M}}{C_2} <  + \infty  $ holds for any fixed $ M>0 $.
\end{theorem}
\begin{remark}
As previously mentioned, one may further consider the case of dependent stochastic processes, which can better capture nonlinear phenomena.
\end{remark}

It is well known that in the general case, the rate in the unweighted central limit theorem has the optimal order $ \mathcal{O}^{\#} \left(N^{-1/2}\right) $,  and assuming  higher order moment may not improve it. However, our weighting function \eqref{CLTFUN} can significantly improve the convergence rate to $ \mathcal{O} \left(N^{-1}\right) $. Another interesting aspect is that our result is optimal at least under the methods in \cite{Kla09}. With \eqref{CLTTHE2} and \eqref{CLTGAILV1} in Section \ref{PPR}, according to the weighted Cauchy-Schwarz inequality, we have 
\[\sum\limits_{n = 1}^N {\theta _n^4}  = \frac{{{{\left( {\theta _1^2} \right)}^2}}}{{{1^1}}} +  \cdots  + \frac{{{{\left( {\theta _N^2} \right)}^2}}}{{{1^1}}} \geqslant \frac{{{{\left( {\theta _1^2 +  \cdots  + \theta _N^2} \right)}^2}}}{{{{\left( {1 +  \cdots  + 1} \right)}^1}}} = \frac{1}{N},\]
 which is precisely the upper bound order we have obtained in Theorem \ref{CLTT1}. However, it is worth emphasizing that the weighting coefficients that can accelerate the rate of the weighted central limit theorem to $ \mathcal{O}\left(N^{-1}\right) $ \textit{may not} necessarily result in an exponential acceleration in the deterministic case, for example, when taking $ {\theta _n} = \sin^2 \left( {\pi n/N} \right)/\sum\nolimits_{j = 1}^N {\sin^2 \left( {\pi j/N} \right)}  $ for $ 1 \leqslant n \leqslant N $ \cite{TL25a}.

This weighted approach can also be extended to establish a variety of additional results, which we recommend exploring further in  \cite[Corollary 4]{MM07},  \cite[Theorem 2.1]{GS09},  \cite[Theorem 1.2]{KS12}, and also  \cite[Theorem 1.1]{Gru14}.  From a probabilistic perspective, under the assumption that the fourth moment is finite (but without the additional requirements in Theorem \ref{CLTT1}), the weights that achieve the weighted central limit theorem at a rate of $ \mathcal{O}(N^{-1}) $ are of  almost full normalized Lebesgue measure $ \mathfrak{s}_{N-1} $ on $ S^{N-1} $ (except for a set of small measure), see   \cite[Theorem 1.1]{KS12}; while under the assumption that the fifth moment is finite, by introducing a corrected normal distribution, it can be shown that the weights that achieve the weighted central limit theorem at a rate of $ \mathcal{O}(N^{-3/2}) $ are of  almost full normalized Lebesgue measure $ \mathfrak{s}_{N-1} $ on $ S^{N-1} $, see   \cite[Theorem 5.1.1]{Bob20}. Unfortunately, we do not know whether our weights derived from  \eqref{CLTFUN} or some other variants (e.g., positive on $ [0,1/2] $ and negative on $ [1/2,1] $, see more general cases in \cite{TL25a}) with the $ C_0^\infty \left([0,1]\right)$ property as well as the normalization property belongs to such almost full measure sets. One feasible approach is to verify certain arithmetic properties of the weights, such as those described in \cite[(iii)]{KS12} (or the assumptions within Conjecture \ref{CONJECTURE} below). This approach requires further exploration and is of significant interest. Here, based on \cite[Theorem 1.1]{KS12} in the measure sense, we leave the following conjecture:

\begin{conjecture}\label{CONJECTURE}
	Assume that $ X_1,  \ldots  , X_N $ are independent and identically distributed  centered random variables, satisfying $ {\bf E}{X_1} = 0,{\bf E}X_1^2 = 1,{\bf E}X_1^3 = 0 $ and $ {\bf E}{\left| {{X_1}} \right|^5} <  + \infty  $. Then, there exists some $ C_0^\infty \left( { [ {0,1}  ]} \right) $ weighting function $ \bar w(x) $ with $ \int_0^1 {\bar w\left( x \right)dx}  = 1 $, such that
for $ {{\bar A}_N} = \sum\nolimits_{n = 0}^{N - 1} {\bar w\left( {n/N} \right)} >0 $ and some universal constant $ C_3>0 $, 
\[\mathop {\sup }\limits_{x \in \mathbb{R}} \left| {{\bf P}\left( {\frac{1}{{\sqrt {{{\bar A}_N}} }}\sum\limits_{n = 1}^{N - 1} {\sqrt {\bar w\left( {n/N} \right)} {X_n}}  \leqslant x} \right) - \mathcal{G}\left( x \right)} \right| \leqslant \frac{{{C_3}{\bf E}{{\left| {{X_1}} \right|}^5}}}{{{N^{3/2}}}},\]
provided the corrected normal distribution $ \mathcal{G}(x)$ defined as
\[\mathcal{G}\left( x \right) = \frac{1}{{\sqrt {2\pi } }}\int_{ - \infty }^x {\exp \left( { - {t^2}/2} \right)dt}  - \frac{{{\bf E}{{\left| {{X_1}} \right|}^4} - 3}}{{8\sqrt {2\pi } N}}\left( {{x^3} - 3x} \right)\exp \left( { - {x^2}/2} \right).\]
\end{conjecture}

From the previous result, it appears that under typical random conditions, the weighting function $ w_{p,q}(x) $ does not accelerate  the convergence rate beyond the original findings, instead maintaining  consistency. This might be due to the highly general scenarios considered, making it difficult to utilize its specific asymptotic properties. However, considering  its impressive  performance in dynamical systems, perhaps within the framework of deterministic  dynamical systems (e.g., \cite{Liv96} based on deterministic hyperbolic  systems)  or probabilistic dynamical systems (e.g., \cite{WL25} based on stochastic Hamiltonian systems), by utilizing the intrinsic properties of these systems, one might achieve results with faster convergence. This is indeed a question that deserves further investigation.

\section{New contributions and historical comparisons}\label{SECCON}
\setcounter{footnote}{0}
\renewcommand{\thefootnote}{\fnsymbol{footnote}}
Because this paper contains both a survey component and a large number of theorems, we choose to add this section to illustrate our new contributions and to compare them with the historical literature, thereby helping readers grasp the novelties.

In the deterministic setting, Theorems \ref{QMC} and \ref{CLT-FOU} are our \textit{main} results, whereas Theorem \ref{LIMIT} is a relatively straightforward result, and Theorem \ref{TH3} serves mainly as an explicit statement rather than a proof.

Theorem \ref{QMC} was conceived with two motivations: first, to summarize the existing results, and second, to establish improved quantitative versions, both of which lend themselves to practical simulations in dynamical systems. More importantly, Theorem \ref{QMC} allows us to investigate, at a theoretical level, how the convergence of weighted Birkhoff averages is shaped by the interplay among the weighting function, nonresonance, and regularity. For these reasons, we prefer to regard \textit{the bulk of} Theorem \ref{QMC} as an \textit{essential} contribution of this paper.
 Specifically, Theorem \ref{QMC} treats the general weighting function $  w_{p,q}(x) $, going beyond the special case $  w_{1,1}(x) $ that has dominated earlier theoretical works \cite{DSSY17,DY18,DM23,TL24a,TL24b,TL25a}. This provides rigorous, quantitative guarantees for the numerical simulations explored in \cite{DSSY17,CCGd24}, which are of considerable practical importance\footnote{It is worth emphasizing that the theoretical analysis becomes far more intricate than the $ p = q = 1 $ case, and the techniques previously employed are no longer effective.}. A detailed comparison with the historical literature will be given in the first paragraph of Section \ref{SECPT12}; here we confine ourselves to a brief summary. Indeed, the polynomial type estimate (I)--(1) is an immediate consequence of the techniques developed in \cite{DSSY17,DY18,DM23,TL24a}, requiring no essential alteration, whereas the exponential type estimate (I)--(3) was established in \cite{TL25b}. Hence, neither forms part of the original contributions of this paper. The exponential type results (I)--(2), (I)--(4), (II)--(5), (II)--(6), and (III)--(7) are \textit{new}. We emphasize, however, that for the special case $ p=q=1 $ the statement corresponding to (I)--(2) was obtained in \cite{TL25b}; therefore we prefer not to regard our treatment of general $ p,q>0 $ as an essential innovation. For the special case $ p=q=1 $, a qualitative version of (II)--(5) was established in \cite{TL24b} (also in a multiple sense); here we improve it \textit{quantitatively} for general $ p,q>0 $, a highly nontrivial yet fundamental challenge that constitutes an \textit{essential} novelty. Similarly, for $ p=q=1 $, a weaker, qualitative form of (II)--(6) was given in \cite{TL24a,TL24b}. Finally, (III)--(7) represents another \textit{essential} novelty: we provide the \textit{first} quantitative exponential convergence result for the purely periodic case. To further illustrate the contrast between our exponential results and earlier ones, consider the most commonly used quasi-periodic case (I)--(2).  When $ p=q=1 $, \cite{TL25b} proves that the weighted Birkhoff average converges at the rate $ \mathcal{O}\left(\exp(-N^{\zeta})\right) $ for any fixed $0< \zeta<(d+2)^{-1} $.  Our work shows that for general $ p,q >0$ the exponent $ \zeta $ depends \textit{explicitly} on these parameters via any fixed $0< \zeta<(d+1+1/\min\{p,q\})^{-1}$, a refinement of both theoretical and practical importance. Similarly, in the most important almost-periodic setting (II)--(6), the qualitative result of \cite{TL24b} guarantees the weighted Birkhoff average to converge at the rate  $  \mathcal{O}\left(\exp\left(-(\log N)^\zeta\right)\right) $ for $ p=q=1 $, yet $ \zeta $ is \textit{only} known to be positive, with \textit{no} explicit connection to the underlying infinite-dimensional spatial structure.  Our new result (II)--(6) makes this dependence transparent as $ 2 \leqslant {\zeta} < 1 + \eta $\footnote{Here $ 2 \leqslant \eta \in \mathbb{N}^+ $ characterizes the properties of the  infinite-dimensional spatial structure.} and is therefore of \textit{significance} in  infinite-dimensional ergodic theory.

Theorem \ref{TH3}, although not regarded as the essential contribution of the present paper---we defer its complete proof to a companion work---is nevertheless a \textit{novel} result. Weighted Birkhoff averages were originally introduced to accelerate the convergence of the classical Birkhoff average, both theoretically and in practical simulations. In all preceding studies (as cited in the Introduction), the weighting functions appearing in Theorem \ref{TH3} and Remark \ref{RE13} have \textit{never} been considered, yet they possess markedly \textit{stronger} theoretical acceleration properties than those currently in use. This reveals the profound link between the asymptotics of the weighting function and the convergence rate of the weighted Birkhoff average, thereby lending the result additional importance.

Theorem \ref{CLT-FOU} is entirely \textit{new} and constitutes one of the \textit{essential} contributions of this paper.  It is worth stressing that this result is intimately tied to the paper's overarching focus on convergence. Weighted Birkhoff averages have proved remarkably effective in concrete computations---notably for calculations of Fourier coefficients, as surveyed in Section \ref{SECphy}---yet \textit{no} rigorous theory has hitherto guaranteed their correctness, let alone established an exponential convergence result.  Moreover, existing applications are confined to the quasi-periodic setting; the almost-periodic case has remained \textit{untouched}.   Motivated by these gaps, we establish an almost complete, quantitative theorem that is simultaneously amenable to practical use.  A key innovation is the introduction of the notion ``\textit{effective order}'', which quantifies, for the \textit{first} time, the range of modes whose Fourier coefficients can be calculated with exponential convergence via weighted Birkhoff averages. This aspect has received virtually \textit{no} prior consideration. Finally, we emphasize that Theorem \ref{CLT-FOU} is \textit{not} a mere corollary of Theorem  \ref{QMC}.  In fact, parts of the proof of the latter are subsumed within the more intricate argument required for the former, which we have explained in Remark \ref{RE14}.

The preceding discussion concerns the novel contributions and comparisons in the deterministic setting.  For the probabilistic case, we regard Theorems \ref{SJT1}, \ref{LIMT1}, \ref{LIMT2} and \ref{CLTT1} less as essential advances than as careful applications of historical results. To maintain objectivity, we have repeatedly acknowledged this point throughout the paper.  Accordingly, these probabilistic statements are not positioned as the paper's main results, even though they appear \textit{new} within the specific framework of weighted Birkhoff averages. Our motivation for establishing them is twofold: they allow us to investigate, from a different perspective, the role of weighted Birkhoff averages in the weighted laws of large numbers, the weighted strong laws of large numbers, and the weighted central limit theorem; and they provide slower convergence rates that stand in sharp contrast to the deterministic case. We believe this perspective opens numerous directions for further investigation.

\section{Proof of the deterministic results}\label{SECPD}
\subsection{Proof of Theorem \ref{LIMIT}}\label{PROOFPRO1}
\begin{proof}
	The proof presented here, though very simple, clearly illustrates the practicality of the weighted Birkhoff average in the sense of limit-preserving. In what follows, we will focus on the discrete case, as the continuous case is essentially the same.
	
	By Abel's summation formula, we have
	\begin{align}
		\frac{1}{{{A^{p,q}_N}}}\sum\limits_{n = 1}^{N } {w_{p,q}\left( {n/N} \right){a_n}}  		& = \frac{1}{{{A^{p,q}_N}}}\left( {w_{p,q}\left( 1 \right)\sum\limits_{j = 1}^N {{a_j}}  + \sum\limits_{n = 1}^{N - 1} {\left( {w_{p,q}\left( {n/N} \right) - w_{p,q}\left( {(n + 1)/N} \right)} \right)\sum\limits_{j = 1}^n {{a_j}} } } \right)\notag \\
		\label{OPTAPP1}	& = \frac{1}{{{A^{p,q}_N}}}\sum\limits_{n = 1}^{N - 1} {n\left( {w_{p,q}\left( {n/N} \right) - w_{p,q}\left( (n+1)/N \right)} \right)\frac{1}{n}\sum\limits_{j = 1}^n {{a_j}} } \\
		& = \frac{N}{{{A^{p,q}_N}}}\sum\limits_{n = 1}^{N - 1} {(n/N)\left( {w_{p,q}\left( {n/N} \right) - w_{p,q}\left( (n + 1)/N \right)} \right)a} \notag \\
		& \;\;\;\;+ \frac{N}{{{A^{p,q}_N}}}\sum\limits_{n = 1}^{N - 1} {(n/N)\left( {w_{p,q}\left( {n/N} \right) - w_{p,q}\left( (n + 1)/N \right)} \right)\left( {\frac{1}{n}\sum\limits_{j = 1}^n {{a_j}}  - a} \right)} \notag \\
		\label{Lim1}: &= {\mathscr{J}_1} + {\mathscr{J}_2}.
	\end{align}
	For the leading part $ \mathscr{S}_1 $, we arrive at
	\begin{align}
		{\mathscr{J}_1}& =  - \frac{N}{{{A^{p,q}_N}}} \cdot \left( {\frac{1}{N}\sum\limits_{n = 1}^{N - 1} {(n/N)w_{p,q}'\left( {{\xi _{n,N}}} \right)} } \right) \cdot a\notag \to  - {\left( {\int_0^1 {w_{p,q}\left( x \right)dx} } \right)^{ - 1}} \cdot \left( {\int_0^1 {xw_{p,q}'\left( x \right)dx} } \right) \cdot a\notag \\
		\label{Lim2}	& =  - {\left( {\int_0^1 {w_{p,q}\left( x \right)dx} } \right)^{ - 1}} \cdot \left( {xw_{p,q}\left( x \right)|_0^1 - \int_0^1 {w_{p,q}\left( x \right)dx} } \right) \cdot a = a,
	\end{align}
	where $ \xi_{n,N} \in (n/N,(n+1)/N) $ for $  1\leqslant n \leqslant N-1 $. As for the remainder $ \mathscr{S}_2 $, we get
	\begin{align}
		{\mathscr{J}_2}  &= \frac{N}{{{A^{p,q}_N}}}\sum\limits_{n = 1}^{N - 1} {(n/N)\left( {w_{p,q}\left( {n/N} \right) - w_{p,q}\left((n + 1)/N \right)} \right)\left( {\frac{1}{n}\sum\limits_{j = 1}^n {{a_j}}  - a} \right)} \notag \\
			& = \frac{N}{{{A^{p,q}_N}}}\sum\limits_{n = 1}^M { + \sum\limits_{n = M + 1}^{N - 1} {(n/N)\left( {w_{p,q}\left( n/N \right) - w_{p,q}\left( (n + 1)/N \right)} \right)\left( {\frac{1}{n}\sum\limits_{j = 1}^n {{a_j}}  - a} \right)} } \notag \\
	\label{Lim3}	:& = \mathscr{J}_2^{\left( 1 \right)} + \mathscr{J}_2^{\left( 2 \right)}.
	\end{align}
	Note that for arbitrarily given $ \varepsilon>0 $, there exists $ M \in {\mathbb{N}^ + } $ such that $ {\| {{n^{ - 1}}\sum\nolimits_{j = 1}^n {{a_j}}  - a} \|_\mathscr{B}} < \varepsilon  $ holds for $ n>M $,	and this leads to $ { \| {{n^{ - 1}}\sum\nolimits_{j = 1}^n {{a_j}}  - a} \|_\mathscr{B}} < U $ with  some $ U>0 $ for all $ n \in \mathbb{N}^+ $. Thus,
	\begin{equation}\label{Lim4}
		{ \| {\mathscr{J}_2^{\left( 1 \right)}}  \|_\mathscr{B}} \leqslant \frac{N}{{{A^{p,q}_N}}}\sum\limits_{n = 1}^M {(n/N) \cdot 2{\left\| w \right\|_{{C^0}\left( {0,1} \right)}}}  \cdot U
		\leqslant \frac{N}{{{A^{p,q}_N}}} \cdot \frac{{UM\left( {M + 1} \right){{ \| w_{p,q}  \|}_{{C^0}\left( {0,1} \right)}}}}{N}
		\leqslant \frac{Q}{N} \to 0
	\end{equation}
	holds with the control coefficient $ Q \lesssim  {\left( {\int_0^1 {w_{p,q}\left( x \right)dx} } \right)^{ - 1}}UM\left( {M + 1} \right){ \| w_{p,q}  \|_{{C^0}\left( {0,1} \right)}} $,	and
	\begin{align}
		{ \| {\mathscr{J}_2^{\left( 2 \right)}}  \|_\mathscr{B}} &\leqslant \sum\limits_{n = M + 1}^{N - 1} {(n/N)\left| {w_{p,q}\left( n/N \right) - w_{p,q}\left( (n + 1)/N \right)} \right|\varepsilon } \notag \\
		\label{Lim5}&	\leqslant \left( {\frac{1}{N}\sum\limits_{n = 1}^{N - 1} {(n/N)\left| {w_{p,q}'\left( {{\xi _{n,N}}} \right)} \right|} } \right) \cdot \varepsilon \to \left( {\int_0^1 {x\left| {w_{p,q}'\left( x \right)} \right|dx} } \right) \cdot \varepsilon .
	\end{align}
	Finally, combining \eqref{Lim1} to \eqref{Lim5}, we obtain the desired conclusion \eqref{OPTPro2.1}.
	
	Next, we construct an example that satisfies \eqref{OPTPro2.1} but contradicts  \eqref{OPTpro1}. Let us consider the Banach space $ ( {\mathbb{R},|  \cdot  |} ) $. Let $ {b_0}: = 0 $ and $ {n^{ - 1}}\sum\nolimits_{j = 1}^n {{a_j}} : = {b_n} $. This leads to $ {a_n} = n{b_n} - (n - 1){b_{n - 1}} $ for all $ n \in \mathbb{N} $. Now, let us make $ b_n $ sparse, e.g., $ {b_n} = 1 $ for $ n = {10^\iota } $ with $ \iota  \in {\mathbb{N}^ + } $, and $ {b_n} = 0 $ for all $ n $ else. Therefore, 	by recalling \eqref{OPTAPP1}, it is evident to verify that
	\[\mathop {\lim }\limits_{N \to  + \infty } \frac{1}{{{A^{p,q}_N}}}\sum\limits_{n = 0}^{N - 1} {w_{p,q}\left( {n/N} \right){a_n}}  = \mathop {\lim }\limits_{N \to  + \infty } \frac{1}{{{A^{p,q}_N}}}\sum\limits_{n = 0}^{N - 1} {n\left( {w_{p,q}\left( {n/N} \right) - w_{p,q}\left( (n + 1)/N \right)} \right){b_n}}  = 0,\]
	but $ \mathop {\lim }_{N \to  + \infty } {N^{ - 1}}\sum\nolimits_{n = 1}^N {{a_n}}  = \mathop {\lim }_{N \to  + \infty } {b_N} $ does not exist. 
	
Finally, note that in the above arguments, we only use the property that the weighting function $  w_{p,q}(x) $ belongs to $ C_0^\infty([0,1]) $. Therefore, the main results of Theorem \ref{LIMIT} also hold for any weighting function in $ C_0^\infty([0,1]) $ as well, with the normalization constant $ A^{p,q}_N $ adjusted accordingly. This completes the proof of  Theorem \ref{LIMIT}.
\end{proof}

\subsection{Proof of  Theorem \ref{QMC}}\label{SECPT12}
\begin{proof}

  	 The result (I)--(1) was established in \cite[Theorem 3.1]{DSSY17}, \cite[Theorem 1.1]{DY18}, \cite[Theorem 3]{DM23}, and \cite[Theorem 2.1]{TL24a}. The result (I)--(2) was initially obtained in \cite[Corollary 2.1]{TL24a} in a qualitative form, where the exponent $ \zeta_1 > 0 $ in $ \mathcal{O}\left( \exp( -N^{\zeta_1} ) \right) $ is small and not explicitly determined. It was subsequently derived in \cite[Theorem 6.1]{TL25b}, in a qualitative manner. For $ p, q > 0 $, the result (I)--(2) can be derived as a byproduct from the proof of Theorem \ref{CLT-FOU}, see Section \ref{SEC14} for details.	 The result (I)--(3) was proved in \cite[Theorem 6.2]{TL25b}. The proof of the result (I)--(4) is similar to that of (III)--(7), which will be discussed later.  	For $ p = q = 1 $, the result (II)--(5) was initially obtained in \cite[Corollary 3.1]{TL24a} in a multiple version and  in a qualitative form. Similarly, for the quantitative result regarding $ p, q > 0 $, it can also be derived from the proof of Theorem \ref{CLT-FOU}. A weaker version of the result (II)--(6) was presented in \cite[Corollary 3.2]{TL24a} and \cite[Theorem 4.4]{TL24b}. Given that the stronger conclusion can be derived by similar arguments (also combining the approach in Theorem \ref{CLT-FOU}), we omit the proof here for brevity. Although some of the aforementioned references only discuss the discrete case, the proof for the continuous case is relatively more straightforward, as pointed out in \cite{TL24a,TL24b} and others. Consequently, we provide the proof only for the result (III)--(7) in this section.

Consider the $ T $-periodic case with $ T \in \mathbb{N}^+ $, namely the data $ f(0),f(1),\cdots,f(T-1) $ with $ f $ being $ T $-periodic. We  first need a \textit{quantitative} version of the inverse discrete Fourier theorem. As it can be seen later, the order of the trigonometric polynomial adapted to the data is somewhat important, since we have to avoid the small divisor to achieve exponential convergence.

For the even-periodic case, i.e., $ T = 2m $ with $ m \in \mathbb{N}^+ $, let the undetermined trigonometric polynomial be
\begin{equation}\label{EVEN}
	f\left( n \right) = {a_0} + \sum\limits_{k = 1}^m {{a_k}\cos \left( {k\pi n/{m}} \right)}  + \sum\limits_{k = 1}^{m - 1} {{b_k}\sin \left( {k\pi n/{m}} \right)} , \quad 1 \leqslant n \leqslant 2m.
\end{equation}
As to the odd-periodic case, i.e., $ T = 2m+1 $ with $ m \in \mathbb{N}^+ $ (the $ 1 $-periodic case is trivial), let the undetermined trigonometric polynomial be
\begin{equation}\label{ODD}
	f\left( n \right) = {a_0} + \sum\limits_{k =1}^m {\left( {{a_k}\cos \left( {2k\pi n/{{(2m + 1)}}} \right) + {b_k}\sin \left( {2k\pi n/{{(m + 1)}}} \right)} \right)} , \quad 1 \leqslant n \leqslant 2m + 1.
\end{equation}
It can be shown that both \eqref{EVEN}
and \eqref{ODD} have unique coefficients solutions $ a_k,b_k $. Moreover, $ a_0 $ is the average of the periodic data, i.e., $ {a_0} = {T^{ - 1}}\sum\nolimits_{n = 0}^{T - 1} {f\left( n \right)}  $. We only give the proof of the  even-periodic case, and the  odd-periodic case is indeed similar.

Denote 
\[\left\{ \begin{gathered}
	f = {\left( {f\left( 0 \right),f\left( 1 \right), \cdots ,f\left( {2m - 1} \right)} \right)^\top}, \hfill \\
	\mathscr{V} = \left( {{\mathscr{V}_0},{\mathscr{V}_1}, \cdots ,{\mathscr{V}_{2m - 1}}} \right), \hfill \\
	{\mathscr{V}_0} = {\left( {1,1, \cdots ,1} \right)^\top}, \hfill \\
	{\mathscr{V}_{2l - 1}} = {\left( {\cos \left( {\pi  \cdot l \cdot 0/m} \right),\cos \left( {\pi  \cdot l \cdot 1/m} \right), \cdots ,\cos \left( {\pi  \cdot l \cdot \left( {2m - 1} \right)/m} \right)} \right)^\top}, \quad 1 \leqslant l \leqslant m, \hfill \\
	{\mathscr{V}_{2l}} = {\left( {\sin \left( {\pi  \cdot l \cdot 0/m} \right),\sin \left( {\pi  \cdot l \cdot 1/m} \right), \cdots ,\sin \left( {\pi  \cdot l \cdot \left( {2m - 1} \right)/m} \right)} \right)^\top}, \quad 1 \leqslant l \leqslant m - 1, \hfill \\
	\mathscr{C} = {\left( {{a_0},{a_1},{b_1}, \cdots ,a_{m-1},{b_{m-1}},{a_{m}}} \right)^\top}. \hfill \\ 
\end{gathered}  \right.\]
Then the equations in \eqref{EVEN} with all $ 1 \leqslant n \leqslant 2m $  yield $ f = \mathscr{V}\mathscr{C} $.

We now establish the orthogonality of vectors in $ \mathscr{V} $, namely $ \mathscr{V}_l \cdot \mathscr{V}_j  =0$ whenever $ l \ne j $. Set $ {\alpha _1} := \pi \left( {l + j} \right)/m $ and $ {\alpha _2} := \pi \left( {l - j} \right)/m $ for $ 1 \leqslant l \leqslant m $ and $ 1 \leqslant j \leqslant m - 1 $; then  
\begin{align*}
	{\mathscr{V}_{2l - 1}} \cdot {\mathscr{V}_{2j}}&  = \sum\limits_{n = 0}^{2m - 1} {\cos \left( {\pi  \cdot l \cdot n/m} \right)\sin \left( {\pi  \cdot j \cdot n/m} \right)} \\
	& = \frac{1}{2}\left( {\sum\limits_{n = 0}^{2m - 1} {\sin \left( {\pi  \cdot \left( {l + j} \right) \cdot n/m} \right)}  - \sum\limits_{n = 0}^{2m - 1} {\sin \left( {\pi  \cdot \left( {l - j} \right) \cdot n/m} \right)} } \right)\\
	& = \frac{1}{2}\left( {\frac{{\sin \left( {\left( {2m - 1} \right){\alpha _1}/2} \right)}}{{\sin \left( {{\alpha _1}/2} \right)}}\sin \left( {m{\alpha _1}} \right) - \frac{{\sin \left( {\left( {2m - 1} \right){\alpha _2}/2} \right)}}{{\sin \left( {{\alpha _2}/2} \right)}}\sin \left( {m{\alpha _2}} \right)} \right)=0,
\end{align*}
because $ \sin \left( {m{\alpha _1}} \right) = \sin \left( {m{\alpha _2}} \right) = 0 $, with the convention    $ \sin \left( {qx} \right)/\sin x :=q $ when $ x=0 $. For    $ {\mathscr{V}_{2l - 1}}$ and $ {\mathscr{V}_{2j - 1}} $ with $ l \ne j $ and $ l + j \leqslant m + \left( {m - 1} \right) = 2m - 1 < 2m $, we obtain from $ {\alpha _1},{\alpha _2} \in \left( {0,\pi } \right) $ that
\begin{align*}
	{\mathscr{V}_{2l - 1}} \cdot {\mathscr{V}_{2j - 1}} &= \sum\limits_{n = 0}^{2m - 1} {\cos \left( {\pi  \cdot l \cdot n/m} \right)\cos \left( {\pi  \cdot j \cdot n/m} \right)} \\
	& = \frac{1}{2}\left( {\sum\limits_{n = 0}^{2m - 1} {\cos \left( {\pi  \cdot \left( {l + j} \right) \cdot n/m} \right)}  + \sum\limits_{n = 0}^{2m - 1} {\cos \left( {\pi  \cdot \left( {l - j} \right) \cdot n/m} \right)} } \right)\\
	& = \frac{1}{2}\left( {\frac{{\cos \left( {\left( {2m - 1} \right){\alpha _1}/2} \right)}}{{\sin \left( {{\alpha _1}/2} \right)}}\sin \left( {m{\alpha _1}} \right) + \frac{{\cos \left( {\left( {2m - 1} \right){\alpha _2}/2} \right)}}{{\sin \left( {{\alpha _2}/2} \right)}}\sin \left( {m{\alpha _2}} \right)} \right) = 0.
\end{align*}
An analogous computation gives  $ {\mathscr{V}_{2l}} \cdot {\mathscr{V}_{2j}}=0 $ for $ l \ne j $ and $ l + j \leqslant 2m - 2 < 2m $:
\begin{align*}
	{\mathscr{V}_{2l}} \cdot {\mathscr{V}_{2j}} &= \sum\limits_{n = 0}^{2m - 1} {\sin \left( {\pi  \cdot l \cdot n/m} \right)\sin \left( {\pi  \cdot j \cdot n/m} \right)} \\
	& = \frac{1}{2}\left( {\frac{{\cos \left( {\left( {2m - 1} \right){\alpha _2}/2} \right)}}{{\sin \left( {{\alpha _2}/2} \right)}}\sin \left( {m{\alpha _2}} \right) - \frac{{\cos \left( {\left( {2m - 1} \right){\alpha _1}/2} \right)}}{{\sin \left( {{\alpha _1}/2} \right)}}\sin \left( {m{\alpha _1}} \right)} \right) = 0.
\end{align*}
Finally, the identity
\begin{equation}\label{GZV0}
	\sum\limits_{n = 0}^{2m - 1} {{\exp\left(l\pi ni/m\right)}}  = \frac{{1 - {\exp\left(2l\pi i\right)}}}{{1 - {\exp\left({l\pi i/m}\right)}}} = 0, \quad \forall 1 \leqslant l \leqslant m
\end{equation} 
provides the orthogonality of $ \mathscr{V}_0 $ and other vectors.

 Note that the  orthogonality implies the uniqueness of $ \mathscr{C} $, so the $ 2m $-periodic data can be expressed as a linear superposition of  trigonometric polynomial of up to order $ m  $. Finally, by summing up all the row vectors in $ \mathscr{V} $ and using \eqref{GZV0}, we obtain that $ {a_0} = {T^{ - 1}}\sum\nolimits_{n = 0}^{T - 1} {f\left( n \right)}  $.

Next, we show the exponential convergence for the weighted periodic case. Denote $ \mathscr{M}_T = T/2 $ for $ T $ even and $ \mathscr{M}_T = \left( {T - 1} \right)/2 $ for $ T $ odd. Then by the previous conclusion, we can rewrite the $ T $-periodic data as a finite order  trigonometric polynomial:
\begin{equation}\label{ZQF}
	f\left( n \right) = {a_0} + \sum\limits_{1 \leqslant \left| k \right| \leqslant \mathscr{M}_T} {{c_k}{\exp\left(2k\pi ni/T\right)}} , \quad \forall n \in \mathbb{N},
\end{equation}
where $ c_k $ are complex coefficients. Note that for all $ 1 \leqslant \left| k \right| \leqslant \mathscr{M}\left( T \right) $, we have $ \left| {k/T} \right| \leqslant \left| {\mathscr{M}_T/T} \right| \leqslant \left| {\left( {T/2} \right)/T} \right| = 1/2 $. This leads to $ \left| {k/T - n} \right| \geqslant 1/2 $ for all $ n \in \mathbb{Z} $. In other words, no small divisors will appear in the following analysis (for more details, see also \cite[Section 3.2]{TL24b}). Then we can conclude from the Poisson summation formula \cite[Chapter 3]{Gra14}  that 
\begin{align}
	\frac{1}{{{A^{p,q}_N}}}\sum\limits_{n = 0}^{N - 1} {w_{p,q}\left( {n/N} \right){\exp\left(2k\pi ni/T\right)}}  &= \frac{1}{{{A^{p,q}_N}}}\sum\limits_{n =  - \infty }^{ + \infty } {w_{p,q}\left( {n/N} \right){\exp\left(2k\pi ni/T\right)}}\notag  \\
	& = \frac{1}{{{A^{p,q}_N}}}\sum\limits_{n =  - \infty }^{ + \infty } {\int_{ - \infty }^{ + \infty } {w_{p,q}\left( {t/N} \right){\exp\left(2k\pi ti/T\right)} \cdot {\exp\left(2\pi nit\right)}dt} }\notag \\
\label{ZQ}	& = \frac{N}{{{A^{p,q}_N}}}\sum\limits_{n =  - \infty }^{ + \infty } {\int_0^1 {w_{p,q}\left( z \right){\exp\left(2\pi iNz\left( {k/T - n} \right)\right)}dz} } ,
\end{align}
which admits exponential smallness $ \mathcal{O}\left( {\exp (-c N^{{\zeta_5}} )} \right) $ for some universal constant $ c>0 $, where $ {\zeta _5} = {\left( {1 + 1/\min \left\{ {p,q} \right\}} \right)^{ - 1}} $, by using similar approach in \cite{TL25b}. Alternatively, we can utilize the main ideas from the proof of Theorem \ref{CLT-FOU} in Section \ref{SEC14}. The distinction lies in the fact that the $ m_N $  defined in \eqref{fbcs} is replaced with $ {m_N}\sim {N^{1/{\beta _{p,q}}}} $ as $ N \to +\infty $, where $ {\beta _{p,q}}: = 1 + 1/\min \left\{ {p,q} \right\} > 1 $, and the convergence rate in \eqref{S11} is replaced with $ \mathcal{O}\left( {\exp (-c N^{{\zeta_5}} )} \right) $. Finally, in view of \eqref{ZQF}, we obtain the exponential convergence for the weighted $ T $-periodic data to its periodic average as
\begin{align*}
	\left| {\frac{1}{{{A^{p,q}_N}}}\sum\limits_{n = 0}^{N - 1} {w_{p,q}\left( {n/N} \right)f\left( n \right)}  - \frac{1}{T}\sum\limits_{s = 0}^{T - 1} {f\left( s \right)} } \right| &= \left| {\sum\limits_{1 \leqslant \left| k \right| \leqslant \mathscr{M}_T} {{c_k}\left( {\frac{1}{{{A^{p,q}_N}}}\sum\limits_{n = 0}^{N - 1} {w_{p,q}\left( {n/N} \right){\exp\left(2k\pi ni/T\right)}} } \right)} } \right|\\
	& \leqslant \mathop {\max }\limits_{1 \leqslant \left| k \right| \leqslant \mathscr{M}_T} \left| {{c_k}} \right| \cdot \sum\limits_{1 \leqslant \left| k \right| \leqslant \mathscr{M}_T} {\mathcal{O}\left( {\exp (-c {{N^{  \zeta_5 }}} )} \right)} \\
	& = \mathcal{O}\left( {\exp ( -c{{N^{  \zeta_5 }}} )} \right),
\end{align*}
provided implied constants that only depend on $ T $ and  $ {\max _{0 \leqslant n \leqslant T}}\left| {f\left( n \right)} \right| $ (here we regard $ p $ and $ q $ as fixed).
\end{proof}

\subsection{Proof of Theorem \ref{CLT-FOU}}\label{SEC14}
\begin{proof}
	Note that almost all rotations $ \rho \in \mathbb{T}^d$ with $  d \in  \mathbb{N}^+ \cup \{  + \infty \} $ are nonresonant,  therefore  the corresponding toral translation is uniquely ergodic on $ \mathbb{T}^d $ \cite{Koz21}. By utilizing the Birkhoff ergodic theorem and Theorem \ref{LIMIT} with $ (\mathbb{R}^D,\|\cdot\|_{\ell^\infty}) $, we conclude that in the $ \|\cdot\|_{\ell^\infty} $ norm, 
\begin{align}
\widehat K\left( s \right) &= \int_{{\mathbb{T}^d}} {K\left( x \right)\exp \left( { - 2\pi i\left\langle {s,x} \right\rangle } \right)dx} \notag \\
& = \mathop {\lim }\limits_{N \to  + \infty } \frac{1}{N}\sum\limits_{n = 1}^{N - 1} {K\left( {{\theta _n}} \right)\exp \left( { - 2\pi i\left\langle {s,{\theta _n}} \right\rangle } \right)} \quad (\text{recall that $ {\theta _n}: = {\theta _0} + n\rho  $}) \notag \\
& = \exp \left( { - 2\pi i\left\langle {s,{\theta _0}} \right\rangle } \right)\mathop {\lim }\limits_{N \to  + \infty } \frac{1}{{A_N^{p,q}}}\sum\limits_{n = 1}^{N - 1} {{w_{p,q}}\left( {n/N} \right)K\left( {{\theta _0} + n\rho } \right)\exp \left( { - 2\pi in\left\langle {s,\rho } \right\rangle } \right)} \notag \\
\label{KBD}& = \exp \left( { - 2\pi i\left\langle {s,{\theta _0}} \right\rangle } \right)\mathop {\lim }\limits_{N \to  + \infty } \frac{1}{{A_N^{p,q}}}\sum\limits_{n = 1}^{N - 1} {{w_{p,q}}\left( {n/N} \right){p_n}\exp \left( { - 2\pi in\left\langle {s,\rho } \right\rangle } \right)} . 
\end{align}
This simply proves \eqref{FOUSL}. As for the more precise version, namely the convergence rate of \eqref{FOUSL}, we need to discuss it based on the universal nonresonance and the spatial structure. Taking account of the previous discussion, we only need to analyze the convergence rate of
\[ \sum\limits_{s \ne k \in {\mathbb{Z}^d}(\mathbb{Z}_ * ^\infty )} {\widehat K\left( k \right)\exp \left( {2\pi i\left\langle {k - s,{\theta _0}} \right\rangle } \right)\left( {\frac{1}{{A_N^{p,q}}}\sum\limits_{n = 1}^{N - 1} {{w_{p,q}}\left( {n/N} \right)\exp \left( {2\pi in\left\langle {k - s,\rho } \right\rangle } \right)} } \right)}. \]
or instead, the convergence rate of
\begin{equation}\label{SN}
	{\mathscr{E}_N}: = \sum\limits_{s \ne k \in {\mathbb{Z}^d}(\mathbb{Z}_ * ^\infty )} {{{\| {\widehat K\left( k \right)} \|}_{{\ell ^\infty }}}\left| {\frac{1}{{A_N^{p,q}}}\sum\limits_{n = 1}^{N - 1} {{w_{p,q}}\left( {n/N} \right)\exp \left( {2\pi in\left\langle {k - s,\rho } \right\rangle } \right)} } \right|} .
\end{equation}
\noindent  \textbf{(I) The finite-dimensional case with analyticity:} Note that almost all rotations $ \rho \in \mathbb{R}^d $ satisfy the finite-dimensional Diophantine nonresonance given in  Definition \ref{Finite-dimensional Diophantine nonresonance}, i.e., 
\[\left| {\left\langle {k,\rho } \right\rangle  - n} \right| > \gamma \left\| k \right\|_{{\ell ^1}}^{ - \tau }, \quad \gamma>0, \quad \forall 0\ne k \in \mathbb{Z}^d, \quad \forall n \in \mathbb{Z},\]
where  $ \tau> d $ can be arbitrarily fixed. Building upon this, for any $ s\in \mathbb{Z}^d $ with $ \|s\|_{\ell^1} \ll {N^{1/\left( {\tau  + {\beta _{p,q}}} \right)}} $,
where the constant $ {\beta _{p,q}} $ is explicitly given by
\[{\beta _{p,q}}: = 1 + 1/\min \left\{ {p,q} \right\} > 1,\]
 we construct the following truncation set and the residual  set, respectively:
\[{\mathcal{S}_{\rm T}^d}: = \left\{ {k \in {\mathbb{Z}^d}: \quad 1 \leqslant {{\left\| {k - s} \right\|}_{{\ell ^1}}} \leqslant {N^{1/\left( {\tau  + {\beta _{p,q}}} \right)}}} \right\},\]
and
\[{\mathcal{S}_{\rm R}^d}: = \left\{ {k \in {\mathbb{Z}^d}: \quad {{\left\| {k - s} \right\|}_{{\ell ^1}}} > {N^{1/\left( {\tau  + {\beta _{p,q}}} \right)}}} \right\}.\]
It is evident to verify that they form a partition of $ {\mathbb{Z}^d}\backslash \left\{ s \right\} $, namely $ {\mathcal {S}_{\rm T}^d} \cup {\mathcal{S}_{\rm R}^d} = {\mathbb{Z}^d}\backslash \left\{ s \right\} $ and $ {\mathcal{S}_{\rm T}^d} \cap {\mathcal{S}_{\rm R}^d} = \phi  $. Therefore, $ {\mathscr{E}_N} $ can be decomposed as
\begin{align}
	{\mathscr{E}_N} &= \sum\limits_{k \in {\mathcal{S}_{\rm T}^d}} { + \sum\limits_{k \in {\mathcal{S}_{\rm R}^d}} {{{\| {\widehat K\left( k \right)} \|}_{{\ell ^\infty }}}\left| {\frac{1}{{A_N^{p,q}}}\sum\limits_{n = 1}^{N - 1} {{w_{p,q}}\left( {n/N} \right)\exp \left( {2\pi in\left\langle {k - s,\rho } \right\rangle } \right)} } \right|} } \notag \\
	:& = \mathscr{E}_N^{(1)} + \mathscr{E}_N^{(2)}. \label{SNDJ} 
\end{align}
In what follows, we will estimate them using different approaches. Specifically, $ \mathscr{E}_N^{(1)} $ involves the difficulty of small divisors caused by nonresonance, while $ \mathscr{E}_N^{(2)} $ is concerned with the effective order of $ s $.

We first estimate $ \mathscr{E}_N^{(1)} $. Similar to the approach in \eqref{ZQ} using the Poisson summation formula \cite[Chapter 3]{Gra14}, we have
\begin{align}
	&\;\left| {\frac{1}{{A_N^{p,q}}}\sum\limits_{n = 1}^{N - 1} {{w_{p,q}}\left( {n/N} \right)\exp \left( {2\pi in\left\langle {k - s,\rho } \right\rangle } \right)} } \right| \notag \\
	= &\;\frac{N}{{A_N^{p,q}}}\left| {\sum\limits_{n =  - \infty }^{ + \infty } {\int_0^1 {{w_{p,q}}\left( z \right)\exp \left( {2\pi iNz\left( {\left\langle {k - s,\rho } \right\rangle  - n} \right)} \right)dz} } } \right|.\label{dengjia}
\end{align}
Based on the observation that for any fixed $k \in  \mathcal{S}_{\rm R}^d $ and nonresonant $ \rho \in \mathbb{R}^d $, there exists a unique $ n^* \in \mathbb{Z}$ such that $ \left| {\left\langle {k - s,\rho } \right\rangle  - {n^ * }} \right| < 1/2 $ which leads to small divisors, we obtain from \cite[Lemma 4.1]{TL25b} (the asymptotic estimate of $ {\left\| {{D^m}{w_{p,q}}} \right\|_{{L^1}\left( {0,1} \right)}} $ for $ m $ sufficiently large) that
\begin{align}
&\;\left| {\int_0^1 {{w_{p,q}}\left( z \right)\exp \left( {2\pi iNz\left( {\left\langle {k - s,\rho } \right\rangle  - {n^ * }} \right)} \right)dz} } \right|\notag \\
 = &\;\left| {{{\left( {2\pi iN\left( {\left\langle {k - s,\rho } \right\rangle  - {n^ * }} \right)} \right)}^{ - m}}\int_0^1 {{D^m}{w_{p,q}}\left( z \right)\exp \left( {2\pi iNz\left( {\left\langle {k - s,\rho } \right\rangle  - {n^ * }} \right)} \right)dz} } \right|\notag \\
 \label{ZJ} \leqslant&\; {\left\| {{D^m}{w_{p,q}}} \right\|_{{L^1}\left( {0,1} \right)}}{\left| {2\pi N\left( {\left\langle {k - s,\rho } \right\rangle  - {n^ * }} \right)} \right|^{ - m}} \\
\leqslant &\;{{\tilde \lambda }^m}{m^{{\beta _{p,q}}m}}{\left( {2\pi N \cdot \gamma \left\| {k - s} \right\|_{{\ell ^1}}^{ - \tau }} \right)^{ - m}}\notag \\
\leqslant&\; \lambda _ * ^m{m^{{\beta _{p,q}}m}}{N^{ - m{\beta _{p,q}}/\left( {\tau  + {\beta _{p,q}}} \right)}}\notag \\
 :=&\;\exp \left( {g\left( m \right)} \right),\notag
\end{align}
where $ \tilde{\lambda}>1 $ is some universal constant, $ {\lambda _ * }: = \tilde \lambda /2\pi \gamma  > 0 $, and the number of integrations by parts $ 2\leqslant  m\in \mathbb{N^+} $ is to be determined. Our goal is to minimize $ \exp \left( {g\left( m \right)} \right) $ in order to achieve the most accurate estimate possible. Therefore, through monotonicity analysis, we set 
\begin{equation}\label{fbcs}
{m_N} \sim \exp \left( { - 1} \right)\lambda _ * ^{ - 1/{\beta _{p,q}}}{N^{1/\left( {\tau  + {\beta _{p,q}}} \right)}} \to  + \infty , \quad m_N \in \mathbb{N}^+,
\end{equation}
thereby obtaining 
\begin{align}
	\left| {\int_0^1 {{w_{p,q}}\left( z \right)\exp \left( {2\pi iNz\left( {\left\langle {k - s,\rho } \right\rangle  - {n^ * }} \right)} \right)dz} } \right|
	& \leqslant \exp \left( {g\left( {{m_N}} \right)} \right) = {\min _{2 \leqslant m \in {\mathbb{N}^ + }}}\exp \left( {g\left( m \right)} \right)\notag \\
	\label{S11} &\lesssim \exp \left( { - {{\beta _{p,q}}}{m_N}} \right) \lesssim  {\exp \left( - c{N^{1/\left( {\tau  + {\beta _{p,q}}} \right)}}\right)}  ,
\end{align}
provided some universal constant $ c>0 $ (here and henceforth). For other $ {{n^ * } \ne n \in \mathbb{Z}} $, the analysis is relatively straightforward because there are no small divisors, and we can simultaneously control the error using \eqref{S11}. More precisely,  with $ 2\leqslant m_N \in \mathbb{N}^+$ defined in \eqref{fbcs}, we derive from \eqref{ZJ} and \eqref{S11} that
\begin{align}
&\;\sum\limits_{{n^ * } \ne n \in \mathbb{Z}} {\left| {\int_0^1 {{w_{p,q}}\left( z \right)\exp \left( {2\pi iNz\left( {\left\langle {k - s,\rho } \right\rangle  - n} \right)} \right)dz} } \right|} \notag \\
 \leqslant &\;\sum\limits_{{n^ * } \ne n \in \mathbb{Z}} {\frac{{{{\left\| {{D^{{m_N}}}{w_{p,q}}} \right\|}_{{L^1}\left( {0,1} \right)}}}}{{{{\left( {2\pi N\left| {\left\langle {k - s,\rho } \right\rangle  - n} \right|} \right)}^{{m_N}}}}}} \notag \\
  \leqslant &\;\frac{{2{{\left\| {{D^{{m_N}}}{w_{p,q}}} \right\|}_{{L^1}\left( {0,1} \right)}}}}{{{{\left( {\pi N} \right)}^{{m_N}}}}} + \frac{{{{\left\| {{D^{{m_N}}}{w_{p,q}}} \right\|}_{{L^1}\left( {0,1} \right)}}}}{{{{\left( {2\pi N} \right)}^{{m_N}}}}}\sum\limits_{{n^ * },{n^ * } \pm 1 \ne n \in \mathbb{Z}} {\frac{1}{{{{\left| {\left\langle {k - s,\rho } \right\rangle  - n} \right|}^{{m_N}}}}}} \notag \\
 \label{S12} \lesssim &\; {\exp \left( { - c{N^{1/\left( {\tau  + {\beta _{p,q}}} \right)}}} \right)} ,
\end{align}
where we have used 
\[\sum\limits_{{n^ * },{n^ * } \pm 1 \ne n \in \mathbb{Z}} {\frac{1}{{{{\left| {\left\langle {k - s,\rho } \right\rangle  - n} \right|}^{{m_N}}}}}}  \leqslant \sum\limits_{{n^ * },{n^ * } \pm 1 \ne n \in \mathbb{Z}} {\frac{1}{{{{\left| {\left\langle {k - s,\rho } \right\rangle  - n} \right|}^2}}}}  \leqslant 2\sum\limits_{n = 0}^{ + \infty } {\frac{1}{{{{\left( {n + 1/2} \right)}^2}}}}  <  + \infty \]
due to
\[\left| {\left\langle {k - s,\rho } \right\rangle  - n} \right| \geqslant \left| {n - {n^ * }} \right| - \left| {\left\langle {k - s,\rho } \right\rangle  - {n^ * }} \right| > \left| {n - {n^ * }} \right| - 1/2, \quad \forall {n^ * },{n^ * } \pm 1 \ne n \in \mathbb{Z}.\]
Now, substituting \eqref{S11} and \eqref{S12} into \eqref{dengjia} and utilizing the simple facts $ \sum\nolimits_{k \in {\mathcal{S}_{\rm T}}} {{{\| {\widehat K\left( k \right)} \|}_{{\ell ^\infty }}}} \leqslant \sum\nolimits_{k \in \mathbb{Z}^d} {{{\| {\widehat K\left( k \right)} \|}_{{\ell ^\infty }}}}  <  + \infty  $\footnote{It is evident that the analyticity of $ K $ implies the absolute summability of its Fourier coefficients.} and $ N/A_N^{p,q} \sim {( {\int_0^1 {{w_{p,q}}\left( x \right)dx} } )^{ - 1}} = 1 $, we obtain 
\begin{equation}\label{S1}
	\mathscr{E}_N^{(1)} \lesssim \exp \left( { - c{N^{1/\left( {\tau  + {\beta _{p,q}}} \right)}}} \right).
\end{equation}

Next, it remains to estimate $ \mathscr{E}_N^{(2)} $. Let $ \sigma>0 $ denote the analytic radius of $ K(x)$ on $ \mathbb{C}^d $; if $ \sigma=+\infty $, it is sufficient to choose $ \sigma=1 $, for instance. Then, given that $ \|s\|_{\ell^1} \ll {N^{1/\left( {\tau  + {\beta _{p,q}}} \right)}} $, and 
\[{N^{1/\left( {\tau  + {\beta _{p,q}}} \right)}} < {\left\| {k - s} \right\|_{{\ell ^1}}} \leqslant {\left\| k \right\|_{{\ell ^1}}} + {\left\| s \right\|_{{\ell ^1}}}, \quad k \in {\mathcal{S}_{\rm R}^d},\]
 it is evident that for $0< \sigma''<\sigma'<\sigma $, we have
\begin{align}
	\mathscr{E}_N^{(2)} &\leqslant \sum\limits_{k \in \mathcal{S}_{\rm R}^d} {{{\| {\widehat K\left( k \right)} \|}_{{\ell ^\infty }}}}  \leqslant {\left\| K \right\|_{\sigma,d} }\sum\limits_{k \in \mathcal{S}_{\rm R}^d} {\exp \left( { - 2\pi \sigma {{\left\| k \right\|}_{{\ell ^1}}}} \right)} \notag \\
	 & = {\left\| K \right\|_{\sigma,d} }\sum\limits_{k \in \mathcal{S}_{\rm R}^d} {\exp \left( { - 2\pi \sigma '{{\left\| k \right\|}_{{\ell ^1}}}} \right) \cdot \exp \left( { - 2\pi \left( {\sigma  - \sigma '} \right){{\left\| k \right\|}_{{\ell ^1}}}} \right)}\notag \\
	 & \lesssim \exp \left( { - 2\pi \sigma '({N^{1/\left( {\tau  + {\beta _{p,q}}} \right)}} - \|s\|_{\ell^1})} \right) \cdot \sum\limits_{k \in \mathcal{S}_{\rm R}^d} {\exp \left( { - 2\pi \left( {\sigma  - \sigma '} \right){{\left\| k \right\|}_{{\ell ^1}}}} \right)} \notag \\
	 & \leqslant \exp \left( { - 2\pi \sigma '({N^{1/\left( {\tau  + {\beta _{p,q}}} \right)}} -  \|s\|_{\ell^1})} \right) \cdot \sum\limits_{k \in {\mathbb{Z}^d}} {\exp \left( { - 2\pi \left( {\sigma  - \sigma '} \right){{\left\| k \right\|}_{{\ell ^1}}}} \right)}  \notag \\
	& \lesssim \exp \left( { - 2\pi \sigma ''{N^{1/\left( {\tau  + {\beta _{p,q}}} \right)}}} \right).\label{S2}
\end{align}
It is observed that the truncation technique employed here can eliminate the nonresonance between $ \rho $ and $k$ within the residual  set $ {\mathcal{S}_{\rm R}^d} $, thereby avoiding the need for small divisor estimates similar to those in $ \mathscr{E}_N^{(1)} $. This idea, initially appeared in \cite{TL24b},  prevents the increase in error, as small divisors can significantly undermine the controllability brought by analyticity as $\|k\|_{\ell^1}$ increases to $+ \infty $.

Finally, by combining \eqref{S1} and \eqref{S2}, we derive  the desired conclusion in (I) from \eqref{SN} and \eqref{SNDJ}. Specifically, we have 
\[ \mathscr{E}_N  \lesssim {\exp ( - {N^{{\zeta _6}}})} , \quad  N \to +\infty\] for any
\[0 < {\zeta _6} < {\left( {d + 1 + 1/\min \left\{ {p,q} \right\}} \right)^{ - 1}},\]
 and this holds uniformly in $ \|s\|_{\ell^1} \ll {N^{{\zeta _6}}}  $. This is valid 
 given that $ \tau>d $ can be arbitrarily fixed in this case. 
\vspace{3mm}

\noindent  \textbf{(II) The infinite-dimensional case with analyticity:}
 In this case, almost all rotations $ \rho \in \mathbb{R}^\infty $ satisfy the infinite-dimensional Diophantine nonresonance presented in  Definition \ref{Infinite-dimensional Diophantine nonresonance}, i.e., 
\begin{equation}\label{WQFGZ}
	\left| {\left\langle {k,\rho } \right\rangle  - n} \right| > \gamma {\left( {\prod\nolimits_{j \in {\mathbb{N} }} {(1 + {{\left\langle j \right\rangle }^\tau }{{\left| {{k_j}} \right|}^\tau })} } \right)^{ - 1}}, \quad \gamma>0, \quad \forall 0\ne k \in \mathbb{Z}^{\infty}_*, \quad \forall n \in \mathbb{Z},
\end{equation}
where  $ \tau> 1 $ can be arbitrarily chosen. For the constant $ 2 \leqslant \eta \in \mathbb{N}^+ $ introduced in the spatial structure, let us fix some $ \zeta_7 $ such that $ 2 \leqslant \zeta_7<1+\eta $. Now, for any $ s \in \mathbb{Z}_*^\infty $ with $ {\left| s \right|_\eta } \ll {( {\log N} )^{\zeta_7} } $, the modified truncation set and the residual  set are then constructed as:
\[\mathcal{S}_{\rm T}^\infty : = \left\{ {k \in \mathbb{Z}_ * ^\infty : \quad 1 \leqslant {{\left| {k - s} \right|}_\eta } \leqslant {{( {\log N} )}^{\zeta_7} }} \right\},\]
and
\[\mathcal{S}_{\rm R}^\infty : = \left\{ {k \in \mathbb{Z}_ * ^\infty : \quad {{\left| {k - s} \right|}_\eta } > {{( {\log N} )}^{\zeta_7} }} \right\}.\]
It can be observed that they naturally form a partition of $ \mathbb{Z}_ * ^\infty \backslash \left\{ s \right\} $. The  framework that follows is similar to (I), yet the details are far more challenging. We still decompose $ 	{\mathscr{E}_N} $ in \eqref{SN} into a form similar to that in \eqref{SNDJ}, on the truncated set $ \mathcal{S}_{\rm T}^\infty  $ and the residual  set $ \mathcal{S}_{\rm R}^\infty  $, and accordingly modify the notations to $  {\tilde{\mathcal{E}}_N^{(1)}} $ and $  {\tilde{\mathcal{E}}_N^{(2)}} $.

Adopting a similar approach, namely decomposing $ \mathbb{Z} $ into a special $  n^*  $ and $ \mathbb{Z} \backslash  \{n^*\} $, we  obtain that
\begin{align}
 {\tilde{\mathcal{E}}_N^{(1)}}& \lesssim {{\tilde \lambda }^{{{\tilde m}_N}}}\tilde m_N^{{\beta _{p,q}}{{\tilde m}_N}}{N^{ - {{\tilde m}_N}}}{\left( {2\pi \gamma \mathop {\inf }\limits_{{{k \in \mathcal{S}_{\rm T}^\infty } }} \prod\limits_{j \in \mathbb{N}} {{{\left( {1 + {{\left\langle j \right\rangle }^\tau }{{\left| {{k_j} - {s_j}} \right|}^\tau }} \right)}^{ - 1}}} } \right)^{ - {{\tilde m}_N}}}\notag \\
& = \lambda _ * ^{{{\tilde m}_N}}\tilde m_N^{{\beta _{p,q}}{{\tilde m}_N}}{N^{ - {{\tilde m}_N}}}{\left( {\mathop {\sup }\limits_{k \in \mathcal{S}_{\rm T}^\infty} \prod\limits_{j \in \mathbb{N}} {\left( {1 + {{\left\langle j \right\rangle }^\tau }{{\left| {{k_j} - {s_j}} \right|}^\tau }} \right)} } \right)^{{{\tilde m}_N}}},\label{WUSD}
\end{align}
where  the number of integrations by parts $ 2\leqslant  \tilde m_N\in \mathbb{N^+} $ is also to be determined.
and we need to estimate the order $ {\mathop {\sup }\nolimits_{k \in \mathcal{S}_{\rm T}^\infty} \prod\nolimits_{j \in \mathbb{N}} {\left( {1 + {{\left\langle j \right\rangle }^\tau }{{\left| {{k_j} - {s_j}} \right|}^\tau }} \right)} } $ in \eqref{WUSD} that  arises from small divisors. For any fixed $ k \in \mathcal{S}_{\rm T}^\infty  $, denote by $ \mathscr{N} $ the number of nonzero components of $ k-s $, and let $ {\left\{ {{k_{{j_\iota }}} - {s_{{j_\iota }}}} \right\}_{1 \leqslant \iota  \leqslant \mathscr{N}}} $ be the corresponding components. Then it follows that
\begin{align}
	{\mathscr{N}^{1 + \eta }} &\lesssim \int_1^\mathscr{N} {{x^\eta }dx}  \lesssim \sum\limits_{\iota  = 1}^\mathscr{N} {{\iota ^\eta }}  \leqslant \sum\limits_{\iota  = 1}^\mathscr{N} {{{\left\langle {{j_\iota }} \right\rangle }^\eta }}  \leqslant \sum\limits_{\iota  = 1}^\mathscr{N} {{{\left\langle {{j_\iota }} \right\rangle }^\eta }\left| {{k_{{j_\iota }}} - {s_{{j_\iota }}}} \right|} \notag \\
& = \sum\limits_{j \in \mathbb{N}} {{{\left\langle j \right\rangle }^\eta }\left| {{k_j} - {s_j}} \right|}  = {\left| {k - s} \right|_\eta } \leqslant {(\log N)^{\zeta_7} }.\label{WU2}
\end{align}
In addition to this, we also have
\begin{align}
{\left\langle {{j_\iota }} \right\rangle ^\tau }{\left| {{k_{{j_\iota }}} - {s_{{j_\iota }}}} \right|^\tau }& = {\left( {\left\langle {{j_\iota }} \right\rangle \left| {{k_{{j_\iota }}} - {s_{{j_\iota }}}} \right|} \right)^\tau } \leqslant {\left( {{{\left\langle {{j_\iota }} \right\rangle }^\eta }\left| {{k_{{j_\iota }}} - {s_{{j_\iota }}}} \right|} \right)^\tau }\notag \\
& \leqslant {\left( {\sum\limits_{j \in \mathbb{N}} {{{\left\langle j \right\rangle }^\eta }\left| {{k_j} - {s_j}} \right|} } \right)^\tau } = \left| {k - s} \right|_\eta ^\tau  \leqslant {(\log N)^{\zeta_7 \tau }},\label{WU3}
\end{align}
because $ \left\langle {{j_\iota }} \right\rangle  = \max \left\{ {1,\left| {{j_\iota }} \right|} \right\} \geqslant 1 $.
Combining \eqref{WU2} and \eqref{WU3} yields that
\begin{align}
	\mathop {\sup }\limits_{k \in \mathcal{S}_{\rm T}^\infty } \prod\limits_{j \in \mathbb{N}} {\left( {1 + {{\left\langle j \right\rangle }^\tau }{{\left| {{k_j} - {s_j}} \right|}^\tau }} \right)}  &= \mathop {\sup }\limits_{k \in \mathcal{S}_{\rm T}^\infty } \exp \left( {\sum\limits_{j \in \mathbb{N}} {\log \left( {1 + {{\left\langle j \right\rangle }^\tau }{{\left| {{k_j} - {s_j}} \right|}^\tau }} \right)} } \right)\notag \\
& = \exp \left( {\mathop {\sup }\limits_{k \in \mathcal{S}_{\rm T}^\infty } \sum\limits_{\iota  = 1}^\mathscr{N} {\log \left( {1 + {{\left\langle {{j_\iota }} \right\rangle }^\tau }{{\left| {{k_{{j_\iota }}} - {s_{{j_\iota }}}} \right|}^\tau }} \right)} } \right)\notag \\
&\leqslant \exp \left( {\mathop {\sup }\limits_{k \in \mathcal{S}_{\rm T}^\infty } \sum\limits_{\iota  = 1}^\mathscr{N} {\log \left( {1 + {{(\log N)}^{\zeta_7 \tau }}} \right)} } \right)\notag \\
&\lesssim \exp \left( {c\mathop {\sup }\limits_{k \in \mathcal{S}_{\rm T}^\infty } \sum\limits_{\iota  = 1}^\mathscr{N} {\log \log N} } \right)\notag \\
	&= \exp \left( {c\mathop {\sup }\limits_{k \in\mathcal{S}_{\rm T}^\infty} \mathscr{N}\log \log N} \right)\notag \\
	& \lesssim \exp \left( {c{{(\log N)}^{\frac{\zeta_7 }{{1 + \eta }}}}\log \log N} \right).\label{WU4} 
\end{align}
Substituting \eqref{WU4} into \eqref{WUSD} and taking 
 \[{{\tilde m}_N} \sim \exp \left( { - 1} \right){N^{1/{\beta _{p,q}}}}\exp \left( {\beta _{p,q}^{ - 1}c{{(\log N)}^{\frac{{\zeta_7} }{{1 + \eta }}}}\log \log N} \right), \quad {\tilde m}_N \in \mathbb{N}^+,\]
it is not too difficult to verify that
\[{{\tilde m}_N} \gtrsim {(\log N)^{\zeta_7} }, \quad {N^a} \gtrsim {(\log N)^{\zeta_7} }\quad  \text{for any $ a>0 $}\]
due to $ 2 \leqslant {\zeta _7} < 1 + \eta  $, as previously required, thus obtaining 
\begin{equation}\label{WU5}
	{\tilde{\mathcal{E}}_N^{(1)}} \lesssim \exp \left( { - {\beta _{p,q}}{{\tilde m}_N}} \right) \lesssim \exp \left( { - c{{(\log N)}^{\zeta_7} }} \right).
\end{equation}

It remains to estimate $ 	{\tilde{\mathcal{E}}_N^{(2)}} $.  We first claim that for any $ b>0 $,
\begin{equation}\label{WU8}
	\sum\limits_{k \in \mathbb{Z}_ * ^\infty } {\exp \left( { - b{{\left| k \right|}_\eta }} \right)}  <  + \infty .
\end{equation}
Once this is proved,  based on this claim  and following the same approach  utilized in \eqref{S2}, we have that for any $0< \sigma''<\sigma'<\sigma $,
\begin{align}
{\tilde{\mathcal{E}}_N^{(2)}}& \lesssim \sum\limits_{k \in \mathcal{S}_{\rm R}^\infty } {{{\| {\widehat K\left( k \right)} \|}_{{\ell ^\infty }}}}  \leqslant {\left\| K \right\|_{\sigma ,\infty }}\sum\limits_{k \in \mathcal{S}_{\rm R}^\infty } {\exp \left( { - 2\pi \sigma {{\left| k \right|}_\eta }} \right)} \notag \\
& = {\left\| K \right\|_{\sigma ,\infty }}\sum\limits_{k \in \mathcal{S}_{\rm R}^\infty } {\exp \left( { - 2\pi \sigma '{{\left| k \right|}_\eta }} \right) \cdot \exp \left( { - 2\pi \left( {\sigma  - \sigma '} \right){{\left| k \right|}_\eta }} \right)} \notag \\
\label{WUQ7}& \leqslant {\left\| K \right\|_{\sigma ,\infty }}\exp \left( { - 2\pi \sigma '\left( {{{(\log N)}^{\zeta_7} } - {{\left| s \right|}_\eta }} \right)} \right)\sum\limits_{k \in \mathbb{Z}_ * ^\infty } {\exp \left( { - 2\pi \left( {\sigma  - \sigma '} \right){{\left| k \right|}_\eta }} \right)} \\ 
& \lesssim \exp \left( { - 2\pi \sigma ''{{(\log N)}^{\zeta_7} }} \right). \label{WUQ6}
\end{align}
Here, \eqref{WUQ7} and \eqref{WUQ6} use the fact that for any $ k \in \mathcal{S}_{\rm R}^\infty  $ and  $ {\left| s \right|_\eta } \ll {( {\log N} )^{\zeta_7} } $ with $ N $ sufficiently large,
\begin{align}
{\left| k \right|_\eta } &= \sum\limits_{j \in \mathbb{N}} {{{\left\langle j \right\rangle }^\eta }\left| {{k_j}} \right|}  \geqslant \sum\limits_{j \in \mathbb{N}} {{{\left\langle j \right\rangle }^\eta }\left( {\left| {{k_j} - {s_j}} \right| - \left| {{s_j}} \right|} \right)} \notag \\
& = {\left| {k - s} \right|_\eta } - {\left| s \right|_\eta } > {(\log N)^{\zeta_7} } - {\left| s \right|_\eta } \geqslant \sigma ''{(\log N)^{\zeta_7} }/\sigma '. \notag
\end{align}
Next, we need to prove the previous claim \eqref{WU8}. In this case, distinct from the  finite-dimensional setting, the infinite-dimensional spatial structure compels us to derive cardinality estimates on the infinite-dimensional lattice, as initially observed in \cite{TL24a,TL24b}. Recall that  \cite[Lemma 8.4]{TL24b} gives 
	\begin{equation}\notag 
	\sum\limits_{ k \in \mathbb{Z}_ * ^\infty ,{{\left| k \right|}_\eta } = \nu } 1  := \# \left\{ {k: \quad  k \in \mathbb{Z}_ * ^\infty, \;{{\left| k \right|}_\eta } = \nu  \in {\mathbb{N}^ + }} \right\} \lesssim {\nu ^{{\nu ^{1/\eta }}}}, \quad  \nu \to +\infty.
\end{equation}
Applying this, it follows from $ \eta \geqslant 2 $ that for any $ b>0 $,
\begin{align*}
\sum\limits_{k \in \mathbb{Z}_ * ^\infty } {\exp \left( { - b{{\left| k \right|}_\eta }} \right)}  &= \sum\limits_{\nu  = 1}^\infty  {\sum\limits_{k \in \mathbb{Z}_ * ^\infty ,{{\left| k \right|}_\eta } = \nu } {\exp \left( { - b{{\left| k \right|}_\eta }} \right)} }  = \sum\limits_{\nu  = 1}^\infty  {\left( {\sum\limits_{k \in \mathbb{Z}_ * ^\infty ,{{\left| k \right|}_\eta } = \nu } 1 } \right)\exp \left( { - b\nu } \right)} \\
& \lesssim \sum\limits_{\nu  = 1}^\infty  {{\nu ^{{\nu ^{1/\eta }}}}\exp \left( { - b\nu } \right)}  = \sum\limits_{\nu  = 1}^\infty  {\exp \left( { - b\nu  + {\nu ^{1/\eta }}\log \nu } \right)} \\
& \leqslant \sum\limits_{\nu  = 1}^\infty  {\exp \left( { - b\nu  + {\nu ^{1/2}}\log \nu } \right)}  \lesssim \sum\limits_{\nu  = 1}^\infty  {\exp \left( { - b\nu /2} \right)}  <  + \infty ,
\end{align*}
as promised.

Finally, combining \eqref{WU5} and \eqref{WUQ6} and recalling \eqref{SN}, we establish the desired universal exponential convergence rate in the infinite-dimensional case via analyticity as
\[	{\mathscr{E}_N} \lesssim {\tilde{\mathcal{E}}_N^{(1)}}+{\tilde{\mathcal{E}}_N^{(2)}} \lesssim {\exp ( { - {{(\log N)}^{\zeta_7} }} )}, \quad  N \to +\infty\]
 for any $ 2 \leqslant \zeta_7<1+\eta $, and this  holds uniformly in $ {\left| s \right|_\eta } \ll {( {\log N} )^{\zeta_7} } $, as $ \zeta_7 $ can be arbitrarily fixed. 
 
 As a remark, the universal infinite-dimensional nonresonance beyond \eqref{WQFGZ} was also constructed in \cite{TL24b}, and the exponential convergence was proved via analyticity (i.e., two distinct approaches were presented there). The method in this section can also be employed to enhance the other approach in \cite{TL24b}.
\vspace{3mm}

The proof of Theorem \ref{CLT-FOU} is now completed.
\end{proof}

\section{Proof of the probabilistic results}\label{PPR}
\subsection{Proof of Theorem \ref{LIMT1}}
\begin{proof}
	Note that $ w_{p,q}(0)=w_{p,q}(1)=0$. Then with $ {\theta _n} := \sqrt {w_{p,q}\left( {n/N} \right)/{A^{p,q}_N}}  $ and the definition of $ A^{p,q}_N $, it is evident that 
	\begin{equation}\label{CLTTHE2}
		\sum\limits_{n = 1}^N {\theta _n^2}  = \sum\limits_{n = 1}^N {w_{p,q}\left( {n/N} \right)/{A^{p,q}_N}}  = \sum\limits_{n = 1}^N { w_{p,q}\left( {n/N} \right)} \Big/\sum\limits_{n = 0}^{N - 1} {w_{p,q}\left( {n/N} \right)}  = 1,  \quad \forall N\geqslant 3.
	\end{equation}

We first consider the case of $ \mu \geqslant 2 $.	By applying \eqref{CLTTHE2} and the previous assumptions, we obtain from \cite[Theorem 3]{Thr87}  the following Marcinkiewcz-Zygmund type strong law of large numbers:
	\[{\bf{P}}\left( {\mathop {\lim }\limits_{N \to  + \infty } \frac{1}{{{N^{1/\mu}}}}\sum\limits_{n = 1}^N {\theta_n {X_n} }  = 0} \right) = 1.\]
	Now, with the  normalization property of the weighting function $ w_{p,q}(x) $, we have
	\begin{equation}\label{NOR}
		{A^{p,q}_N} = \sum\limits_{n = 0}^{N - 1} {{w_{p,q}}\left( {n/N} \right)}  = N\sum\limits_{n = 0}^{N - 1} {{N^{ - 1}}{w_{p,q}}\left( {n/N} \right)}  \sim N\int_0^1 {{w_{p,q}}\left( x \right)dx}  = N, \quad N \to  + \infty .
	\end{equation}
	Therefore, by taking $ {\sigma} :=  1/\mu+ 1/2 $, we prove the desired conclusion for $ \mu>2 $ as
		\[{\bf{P}}\left( {\mathop {\lim }\limits_{N \to  + \infty }  {\frac{1}{{{N^{ \sigma }}}}\sum\limits_{n = 1}^N {\sqrt {{w_{p,q}}\left( {n/N} \right)} {X_n}} } =0} \right) = 1.\]

	It remains to consider the case of $ 0<\mu \leqslant 2 $.	Note that for $ {b_{N,n}}: = \sqrt {{w_{p,q}}\left( {n/N} \right)} /N $, we have 
\[\mathop {\sup }\limits_{N \geqslant 1} \sum\limits_{n = 1}^N {\left| {{b_{N,n}}} \right|}  = \mathop {\sup }\limits_{N \geqslant 1} \frac{1}{N}\sum\limits_{n = 1}^N {\sqrt {{w_{p,q}}\left( {n/N} \right)} \sim} \int_0^1 {w_{p,q}^{1/2}\left( x \right)dx}  <  + \infty .\]
Therefore, with $ \sigma : = 1/\mu  + 1 $, \cite[Lemma 1]{Thr87} provides  
\[{\bf{P}}\left( {\mathop {\lim }\limits_{N \to  + \infty } \frac{1}{{{N^{\sigma}}}}\sum\limits_{n = 1}^N {\sqrt {{w_{p,q}}\left( {n/N} \right)} {X_n}}  = 0} \right) = {\bf{P}}\left( {\mathop {\lim }\limits_{N \to  + \infty } \frac{1}{{{N^{1/\mu}}}}\sum\limits_{n = 1}^N {\left| {{b_{N,n}}} \right|{X_n}}  = 0} \right) = 1,\]
as desired.
\end{proof}

\subsection{Proof of Theorem \ref{LIMT2}}
\begin{proof} 
	Denote $ \mu': =2\mu/(\mu-2)>0 $ with $ \mu>2 $, hence we have $ 1/2 = 1/\mu' + 1/\mu $.  For $ {a_{N,n}}: = \sqrt {{w_{p,q}}\left( {n/N} \right)}  $, we have for any $ \omega  > 0 $ that 
\begin{align*}
{A_\omega }: &= \mathop { {\lim\sup } }\limits_{N \to  + \infty } {\left( {\frac{1}{N}\sum\limits_{n = 1}^N {a_{N,n}^\omega } } \right)^{1/\omega }} = \mathop { {\lim \sup} }\limits_{N \to  + \infty } {\left( {\frac{1}{N}\sum\limits_{n = 1}^N {{{\left( {\sqrt {{w_{p,q}}\left( {n/N} \right)} } \right)}^\omega }} } \right)^{1/\omega }}\\
&= {\left( {\int_0^1 {w_{p,q}^{\omega /2}\left( x \right)dx} } \right)^{1/\omega }} <  + \infty ,
\end{align*}
as the weighting function $ w_{p,q}(x) $ does not admit any singularity.	Therefore, by utilizing	\cite[Theorem 1.4]{CG07} with $ \omega  = \mu' $ (i.e., $ A_{\mu'}<+\infty $) and the previous assumptions, we obtain 
\begin{align*}
&{\bf{P}}\left( {\mathop { {\lim\sup } }\limits_{N \to +\infty  } \left| {\frac{1}{{\sqrt {N\log N} }}\sum\limits_{n = 1}^N {\sqrt {{w_{p,q}}\left( {n/N} \right)} {X_n}} } \right| \leqslant \sqrt {2A_2^2{\bf{E}}{X_1^2}} } \right)\\
 =\; & {\bf{P}}\left( {\mathop { {\lim \sup} }\limits_{N \to +\infty } \left| {\frac{1}{{\sqrt {N\log N} }}\sum\limits_{n = 1}^N {{a_{N,n}}{X_n}} } \right| \leqslant \sqrt {2A_2^2{\bf{E}}{X_1^2}} } \right)= 1,
\end{align*}
	which proves the result with $ C_1:=\sqrt 2 {A_2} > 0 $.
\end{proof}

\subsection{Proof of Theorem \ref{CLTT1}}
\begin{proof}
Recall from \eqref{CLTTHE2} that for $ {\theta _n} := \sqrt {w_{p,q}\left( {n/N} \right)/{A^{p,q}_N}}  $, we have $ \sum\nolimits_{n = 1}^N {\theta _n^2}  =1 $. Furthermore, by	employing \eqref{NOR}, we can further estimate that
	\begin{align}
		\sum\limits_{n = 1}^N {\theta _n^4} & = \sum\limits_{n = 1}^N {w_{p,q}^2\left( {n/N} \right)/(A^{p,q}_N)^2}  \sim \sum\limits_{n = 1}^N {w_{p,q}^2\left( {n/N} \right)/{N^2}}\notag \\ 
		\label{CLTTHE4}	&= {N^{ - 1}}\sum\limits_{n = 1}^N {{N^{ - 1}}w_{p,q}^2\left( {n/N} \right)} 
		\sim {N^{ - 1}}\int_0^1 {w_{p,q}^2\left( x \right)dx} .
	\end{align}
	Now, by \eqref{CLTTHE2} and the previous assumptions,  applying  \cite[Theorem 1]{Kla09} yields
	\begin{equation}\label{CLTGAILV1}
		\mathop {\sup }\limits_{\alpha ,\beta  \in \mathbb{R},\;\alpha  < \beta } \left| {{\bf{P}}\left( {\alpha  \leqslant \frac{1}{{\sqrt {{A^{p,q}_N}} }}\sum\limits_{n = 1}^N {\sqrt {w_{p,q}\left( {n/N} \right)} {X_n}}  \leqslant \beta } \right) - \frac{1}{{\sqrt {2\pi } }}\int_\alpha ^\beta  {{\exp\left(- t^2/2\right)}dt} } \right| \leqslant \tilde{C}_2\sum\limits_{n = 1}^N {\theta _n^4} ,
	\end{equation}
	provided a universal constant $  \tilde{C}_2>0 $. Finally, substituting \eqref{CLTTHE4} into \eqref{CLTGAILV1}, we obtain the desired conclusion. Here, we mention that the final control coefficient $ C_2>0 $ may depend on $ p,q $, and it may tend to infinity as $ p,q\to +\infty $. However, it can be shown in \eqref{CLTTHE4} that $ {\sup _{0 < p,q \leqslant M}}{C_2} <  + \infty  $ holds for any fixed $ M>0 $, as by Lebesgue's Theorem,
	\[\mathop {\lim }\limits_{\left( {p,q} \right) \to \left( {0,0} \right)} \int_0^1 {w_{p,q}^2\left( x \right)dx}  = {\left( {\int_0^1 {{\exp(-1)}ds} } \right)^{ - 2}}\int_0^1 {{\exp(-2)}dx}  = 1.\]
	This completes the proof of Theorem \ref{CLTT1}.
\end{proof}

\section*{Declarations}

\subsection*{Conflict of interest} On behalf of all authors, the corresponding author states that there is no conflict of interest.

\subsection*{Data availability statements}
Data will be made available on reasonable request.

 \section*{Acknowledgements} 
The authors would like to sincerely thank the anonymous referees for their valuable suggestions and comments, which significantly improved the paper.
 Z. Tong sincerely thanks Prof. James D. Meiss (University of Colorado) for his generous sharing of the latest developments in weighted Birkhoff averages through insightful discussions.  Z. Tong also sincerely thanks Prof. Valery V. Ryzhikov (Moscow State University) for his discussions on the weighted rates of general dynamical systems. Z. Tong sincerely thanks Prof. Alessandra Celletti (Universit\`a Roma Tor Vergata) for her 2024 talk related to weighted Birkhoff averages. Z. Tong is also very grateful to Profs. Ai Hua Fan (Universit\'e de Picardie) and Konstantin M. Khanin (University of Toronto) for the discussions, which have helped to extend some of the ideas from the initial version of this paper. These extensions will be discussed in a separate work.  Z. Tong sincerely thanks Prof. Bo'az B. Klartag (Weizmann Institute of Science) for providing valuable information on the advances in weighted central limit theorems.  Z. Tong sincerely thanks Profs.   Aleksandr G. Kachurovski\u{\i} (Sobolev Institute of Mathematics), Michael Lin (Ben-Gurion University), Ivan V. Podvigin (Sobolev Institute of Mathematics) and Naci Saldi (Bilkent University) for their invaluable assistance and discussions.   Z. Tong  was supported by the China Postdoctoral Science Foundation (Grant No. 2025M783102). Y. Li was supported in part by the National Natural Science Foundation of China (Grant Nos. 12071175, 12471183 and 12531009).

\end{document}